\documentclass{article}

\usepackage{amsthm}

\newtheoremstyle{new_plain}
  {0.5\topsep}   
  {0.5\topsep}   
  {\itshape}  
  {0em}       
  {\bfseries} 
  {.}         
  {5pt plus 1pt minus 1pt} 
  {}          
\theoremstyle{new_plain}
\newtheorem{theorem}{Theorem}
\newtheorem{lemma}{Lemma}
\newtheorem{assumption}{Assumption}
\newtheorem{corollary}{Corollary}

\newtheoremstyle{new_definition}
  {0.5\topsep}   
  {0.5\topsep}   
  {\normalfont}  
  {0pt}       
  {\bfseries} 
  {.}         
  {5pt plus 1pt minus 1pt} 
  {}          
\theoremstyle{new_definition}
\newtheorem{definition}{Definition}

\newtheoremstyle{new_remark}
  {0.5\topsep}   
  {0.5\topsep}   
  {\normalfont}  
  {1em}       
  {\itshape} 
  {.}         
  {5pt plus 1pt minus 1pt} 
  {}          
\theoremstyle{new_remark}
\newtheorem{remark}{Remark}

\usepackage[dvipsnames]{xcolor}

\usepackage{graphicx} 
\usepackage{epsfig} 
\usepackage{mathptmx} 
\usepackage{mathptmx} 
\usepackage{amsmath} 
\usepackage{amssymb}  
\usepackage{xcolor}
\usepackage[cal=cm,scr=boondox]{mathalpha}

\DeclareRobustCommand{\svdots}{
  \vcenter{%
    \offinterlineskip
    \hbox{.}
    \vskip0.1\normalbaselineskip
    \hbox{.}
    \vskip0.075\normalbaselineskip
    \hbox{.}%
  }%
}

\usepackage[papersize={210mm,297mm}]{geometry} 
\usepackage{amsmath}
\usepackage{graphicx}
\usepackage{hyperref}
\usepackage{enumitem}
\usepackage{amssymb}
\title{On robotic manipulators with time-dependent inertial parameters: From physical consistency to boundedness of the mass matrix}
\author{Tom Kaufmann\footnote{Gratefully acknowledges support by Deutsche Forschungsgemeinschaft in the framework of Research Training Group $2182$.}\hspace{1.25mm} and Johann Reger\footnote{Both authors are with the Control Engineering Group at Technische Universität Ilmenau, P.O. Box $10$ $05$ $65$, D-$98684$ Ilmenau, Germany.}} 

\begin{document}
\maketitle
\section*{Abstract} We generalize the robotics equation describing the dynamics of open kinematic chains by including the effect of time-dependent change of inertial parameters as well as the effects of causative mass-density redistribution, triggered by internal movement of mass-carrying particles relative to their body-fixed frames. Time dependency of inertial parameters that results from the sole addition of mass to the robot prominently occurs during the loading of end-effectors---a scenario covered by our model without restriction from the restraint that kinematic parameters of the robot must remain constant. Further, our model also includes internal mass-density redistributions that adhere to this kinematic restraint such as trolleys attached to the robot or the movement of passengers. To accompany the generalized robotics equation with some theoretical infrastructure, we then introduce the concepts of uniform physical consistency and upper boundedness of inertial parameters under which desirable, structural properties regarding the existence of finite, positive uniform bounds of the mass matrix can be shown to carry over to the more involved case of time-dependent inertial parameters. These findings have implications for adaptive control, as they facilitate more realistic testing for robustness against unforeseen time dependencies. Moreover, the results in this paper also provide a pathway to ensuring the desirable existence of finite, positive uniform bounds of the estimated mass matrix under upper bounded, uniformly physically consistent estimation regimes. 

\section{Introduction}
Classical approaches to adaptive control of robotic manipulators such as \cite{craig_adaptive_1987,slotine_adaptive_1987,johansson_adaptive_1990} are concerned with the design of feedback to enable tracking when some parameters $\Theta$ of the manipulator are not known. Usually, adaptive structures include an update law that provides suitable estimates $\hat{\Theta}$ of the unknown parameters leading to tracking by means of a control law $\tau$ that uses these estimates instead of the unknown parameters. A specific type of parameters describing any robotic manipulator that received some detailed attention are the inertial ones, where works like e.g.\ \cite{wu_adaptive_2022,patnaik_adaptive_2023,cho_recursive_2024} building on the pioneering contributions \cite{lee_natural_2018,lee_geometric_2019} enforce physical consistency \cite{traversaro_identification_2016,wensing_linear_2018} of the estimates in addition to facilitate tracking. This is done with the promise that physical consistency of the estimates will yield more refined adaptive control action. Anyhow, before being applicable safely in practice, adaptive schemes should  be evaluated for their robustness against unintended uncertainties, e.g.\ of \cite{craig_adaptive_1987,slotine_adaptive_1987,johansson_adaptive_1990,lee_natural_2018,wu_adaptive_2022} against disturbance $w$ and time-dependent change of supposedly constant unknown parameters, i.e.\ against $\dot{\Theta}(t)\neq 0$ for some $t\geq 0$. This makes it necessary to scrutinize them with respect to a model that describes sufficiently well the effects of these unintended uncertainties on the dynamical behavior of the robotic manipulator. Starting from the classical robotics equation
\begin{align}
\label{ClassicalRoboticEquation}
    M(q,\Theta)\ddot{q} + C(q,\dot{q},\Theta)\dot{q}+G(q,\Theta)=&\tau+w
\end{align}
that describes the dynamics of the robotic manipulator when $\Theta$ is constant, the authors of \cite{pagilla_adaptive_2000} work with the more detailed model
\begin{align}
\label{SimpleGeneralizationOfTheRoboticEquation}
    M(q,\Theta)\ddot{q} + \left( C(q,\Theta)+\mathcal{F}(q,\Theta,\dot{\Theta}) \right)\dot{q}+G(q,\Theta)=&\tau+w
\end{align}
of this behavior when $\Theta$ depends on time. This model can be obtained by application of the Lagrange formalism with respect to the kinetic energy $\mathcal{T}$ and the potential energy $\mathcal{U}$ of the robotic manipulator, taking into account that possibly $\dot{\Theta}(t)\neq 0$ for some $t\geq 0$, but using the same structure of $\mathcal{T}$ and $\mathcal{U}$ as if $\Theta$ were constant. This is a suitable approach when investigating the influence of general parameters on the robotic manipulator. However, since many adaptive approaches are concerned with inertial parameters specifically, the model should also include the effect of the cause for their time dependency, i.e.\ time-dependent change of the mass densities, originating from movements of mass-carrying particles\footnote{Please note that here in this work, the term \textit{mass-carrying particle} refers to an entity without spatial extent in 3-dimensional space that contributes to the value of the mass density of the body built by this particle, evaluated at the location of the particle.} that otherwise would remain standing still when deriving the classical robotics equation~\eqref{ClassicalRoboticEquation}. As we will show, these movements introduce separate terms to the kinetic energy, which are ignored in \cite{pagilla_adaptive_2000}, thereby leading to a more realistic version of the robotics equation compared to~\eqref{SimpleGeneralizationOfTheRoboticEquation} when deriving the dynamical behavior by means of the Lagrange formalism.

\subsection{Contribution and outline}
After introducing the style of notation in Section~\ref{section:Notation}, we highlight in Section~\ref{section:PhysicalConsistency} the intrinsic relationship between mass density and inertial parameters of a rigid body---showcasing the fact that time dependency of the inertial parameters is not possible without time-dependent change of the  mass density---and reiterate the related concept of physical consistency. In the same breath, we remove ambiguity present in prior literature regarding the properties of physically consistent inertial parameters by taking a measure-theoretic approach to the matter. Then, we derive in Section~\ref{section:Modeling} a generalization of the robotics equation that includes all relevant effects associated with time-dependent change of the inertial parameters, answering the need for realistic testing by providing a testbed for robustness of adaptive algorithms concerned with counteracting imprecise knowledge of inertial parameters. We explain assumptions to be made such that the dynamics result in an ODE, thus allowing analysis with similar Lyapunov tools as are typically used when working on the classical robotics equation~\eqref{ClassicalRoboticEquation} and therefore striking a middle ground between complexity and realism. As it turns out, this approach describes the effect of time-dependent change of the mass densities of end-effectors without any restriction from the imposed assumptions, thus enabling  analysis of robustness against loading processes without resorting to a PDE machinery, i.e.\ alleviating computational and analytical cost. In Section~\ref{section:StructuralProperties}, we use the insight from the modeling to reveal structural properties of the generalized robotics equation starting with some inherent properties regarding skew symmetry and decomposition into regressor form. Afterward, the notion of physical consistency of inertial parameters of rigid bodies, i.e.\ constant ones, as established in Section~\ref{section:PhysicalConsistency} serves as fundament for the findings in Section~\ref{section:BoundedMassMatrix} on the boundedness of the mass matrix $M(q,\Theta)$ depending on the physical consistency of the inertial parameters. We arrive at generalizations of statements existing in the literature \cite{ghorbel_uniform_1998} for constant $\Theta$ to the time-dependent case, whereby we introduce the concepts of {uniform} physical consistency and upper boundedness of inertial parameters which we then show to be crucial for the existence of finite, positive {uniform} bounds of the mass matrix. These results are interesting in their own right---regardless whether the inertial parameters of a given robotic manipulator are constant or not---because they pave a way to ensure the existence of finite, positive uniform bounds of the estimated mass matrix $\hat{M}(q,t)={M}(q,\hat{\Theta}(t))$ by means of physically consistent adaptation as proposed in \cite{lee_natural_2018} or subsequent, related work \cite{wu_adaptive_2022,patnaik_adaptive_2023,cho_recursive_2024}. Next, we leverage in Section~\ref{section:OtherComponents} the structural insight from the modeling to derive conditions under which other components in the generalized robotics equation that originate from internal movement of mass-carrying particles are bounded. Specifically, this leads to statements
connecting the velocity and acceleration of such movements as well as the rate of change of the inertial parameters to safe
operation of the robotic manipulator. Finally, we offer concluding remarks as well as an outlook to future work in Section~\ref{section:Conclusion}.

\subsection{Notation}\label{section:Notation}
The identity matrix of dimension $p$ is written as $\mathrm{I}_p$. Let $\mathrm{sym}(p)=\{M\in\mathbb{R}^{p\times p}:M=M^\top\}$ and $\mathrm{skew}(p)=\{M\in\mathbb{R}^{p\times p}:M=-M^\top\}$. For any $A\in\mathrm{sym}(p)$, its trace is denoted by $\mathrm{tr}(A)=\sum_{i=1}^p\lambda_i(A)$ and the extremal eigenvalues are represented by $\lambda_{\min}(A)$ and $\lambda_{\max}(A)$. The maximal singular value of any $B\in\mathbb{R}^{p_1\times p_2}$ is given by $\sigma_{\max}(B)=\sqrt{\lambda_{\max}(B^\top B)}$. We note matrices and vectors as
\begin{align}
(r_k)_{k\in\{1,\ldots,p\}}&=\left({\begin{smallmatrix}
    r_1\\
    \svdots\\
    r_{p}
\end{smallmatrix}}\right),\quad (M_{k_1,k_2})_{(k_1,k_2)\in\{1,\ldots,p_1\}\times\{1,\ldots,p_2\}}=\left({\begin{smallmatrix}
    M_{1,1}&\ldots&M_{1,p_2}\\
    \svdots& &\svdots\\
    M_{p_1,1}&\ldots&M_{p_1,p_2}
\end{smallmatrix}}\right),
\end{align}
respectively. Further, let 
\begin{align}
\mathcal{S}(x)&=\left({\begin{smallmatrix}
    0&-x_3&x_2\\
    x_3&0&-x_1\\
    -x_2&x_1&0
\end{smallmatrix}}\right)\in\mathrm{skew}(3)
\end{align}
with $x=\left(x_k\right)_{k\in\{1,2,3\}}\in\mathbb{R}^3$ such that $\mathcal{S}(x)y=x\times y$ for all $y\in\mathbb{R}^3$.

\section{Revisiting the concept of physical consistency of inertial parameters of rigid bodies}\label{section:PhysicalConsistency}
A rigid body is defined by its mass density $\rho(x)\in\mathbb{R}$ in a body-fixed frame with coordinates $x\in\mathbb{R}^3$. Its ten inertial parameters are the mass $m\in\mathbb{R}$, the first moment of mass \begin{align}
\label{eq:h}
h=\begin{pmatrix}h_1&h_{2}&h_{3}\end{pmatrix}^\top\in\mathbb{R}^3
\end{align}
and the inertia matrix
\begin{align}
\label{eq:I}
    I=&\left(\begin{smallmatrix}
        I_{11}&I_{{12}}&I_{{13}}\\
        I_{{12}}&I_{{22}}&I_{{23}}\\
        I_{{13}}&I_{{23}}&I_{{33}}
    \end{smallmatrix}\right)\in\mathrm{sym}(3)
\end{align}
 for rotation about the origin of the body-fixed frame. These inertial parameters are collected in the vector
\begin{align}
    \label{eq:Phi}
    \Phi^\top = &\begin{pmatrix}m  & h_{{1}}  & h_{{2}}  & h_{{3}}  & I_{{11}}  & I_{{22}}  & I_{{33}}  & I_{{12}}  & I_{{23}}  & I_{{13}}\end{pmatrix}^\top \in \mathbb{R}^{10}
\end{align}
and can be calculated from the mass density $\rho$ as follows:
\begin{align}
\label{MassDistributionToInertialParam}
         m=&\int_{\mathbb{R}^3} \rho(x)\mathrm{d}x,\quad
        h=\int_{\mathbb{R}^3}\ \rho(x)x\mathrm{d}x,\quad
        I=\int_{\mathbb{R}^3} \rho(x)\mathcal{S}(x)^\top \mathcal{S}(x)\mathrm{d}x.
\end{align}
A rigid body is physically meaningful if it is defined by a nonnegative mass density, i.e.\ $\rho(x)\geq 0$ for all $x\in\mathbb{R}^3$, such that~\eqref{MassDistributionToInertialParam} yields $m>0$. The fact that there are inertial parameters $\Phi\in\mathbb{R}^{10}$ with $m>0$ for which~\eqref{MassDistributionToInertialParam} remains unsatisfied for all nonnegative mass densities motivates the notion of physical consistency of inertial parameters presented in \cite{wensing_linear_2018} that calls for $m>0$ as well as realizability through some nonnegative mass density such that~\eqref{MassDistributionToInertialParam} is satisfied. According to (\cite{wensing_linear_2018}, Theorem 3), physical consistency of given inertial parameters $\Phi\in\mathbb{R}^{10}$ with entries as in~\eqref{eq:Phi} is equivalent to positive definiteness of their symmetric $4\times 4$ matrix representation $f(\Phi)$, where the function $f:\mathbb{R}^{10}\to\mathrm{sym}(4)$ is defined as
\begin{align}
\label{DefTrafoInertialParameters}
    f(\Phi)=&\begin{pmatrix}
        \Sigma&h\\h^\top&m
    \end{pmatrix},\quad \Sigma=\frac{1}{2}\mathrm{tr}(I)\mathrm{I}_3-I
\end{align}
with $h \in \mathbb{R}^3$ and $I \in \mathrm{sym}(3)$ structured as per $(\ref{eq:h})$, $(\ref{eq:I})$, respectively. This equivalence is a concise characterization of physical consistency that provides the fundament on which we rest our structural findings on the boundedness of the mass matrix in the case of time-dependent inertial parameters in Section~\ref{section:StructuralProperties}. 

Before we do that, we highlight a gap in the proof that is provided by \cite{wensing_linear_2018} for this equivalence, i.e.\  for the claim that physical consistency of $\Phi$ is exactly characterized by positive definiteness of $f(\Phi)$. This requires dissecting the underlying mechanism that points toward such a result. We then propose a refined formulation of the definition of physical consistency that allows bringing forth a novel measure-theoretic argument rigorously proving the equivalence of physical consistency of $\Phi$ and positive definiteness of $f(\Phi)$. First of all, we reformulate the latter of the two statements: Due to Schur's complement, we have
\begin{align}
\label{fPosDef_IFF_PosMassPosDefSigmaCoM}
f(\Phi)\succ 0\quad\overset{\eqref{DefTrafoInertialParameters}}{\Longleftrightarrow}\quad& m>0\text{ and }\underbrace{\Sigma-m^{-1}h h^\top}_{=\Sigma_\mathrm{CoM}}\succ 0,
\end{align}
where $\Sigma_\mathrm{CoM}$ is identified in (\cite{wensing_linear_2018}, Eq.\ $(20)$) as \textit{density-weighted
covariance} associated with rigid bodies for which the quantities $m\in\mathbb{R}$, $\Sigma\in\mathrm{sym}(3)$, $h\in\mathbb{R}^3$ exist and $m>0$. Further, given $m>0$, we can introduce $I_\mathrm{CoM}=I-m^{-1}\mathcal{S}(h)^\top\mathcal{S}(h)\in\mathrm{sym}(3)$, which, in view of the parallel axis theorem (\cite{wensing_linear_2018}, Eq.\ (8)), is the inertia matrix for rotation of a rigid body around its \textit{center of mass} (CoM) located at $c=m^{-1}h\in\mathbb{R}^3$. Now, under $m>0$ and via exploitation of the geometric fact $\mathcal{S}(h)^\top\mathcal{S}(h)=\mathrm{tr}(hh^\top)\mathrm{I}_3-hh^\top$, we can rewrite
\begin{align}
    \Sigma_\mathrm{CoM}&\underset{\Sigma_\mathrm{CoM},\, \Sigma}{\overset{\text{Def. of}}{=}}\frac{1}{2}\mathrm{tr}(I)\mathrm{I}_3-I-m^{-1}hh^\top\allowdisplaybreaks\\
    &\underset{I_\mathrm{CoM}}{\overset{\text{Def. of}}{=}}\frac{1}{2}\mathrm{tr}\left(I_\mathrm{CoM}+m^{-1}\mathcal{S}(h)^\top\mathcal{S}(h)\right)\mathrm{I}_3-I_\mathrm{CoM}-m^{-1}\left(\mathcal{S}(h)^\top\mathcal{S}(h)+hh^\top\right)\allowdisplaybreaks\\
    &\underset{\text{fact}}{\overset{\text{geometric}}{=}}\frac{1}{2}\mathrm{tr}\left(I_\mathrm{CoM}+m^{-1}\left(\mathrm{tr}(hh^\top)\mathrm{I}_3-hh^\top\right)\right)\mathrm{I}_3-I_\mathrm{CoM}-m^{-1}\mathrm{tr}(hh^\top)\mathrm{I}_3\allowdisplaybreaks\\
    &=\frac{1}{2}\mathrm{tr}\left(I_\mathrm{CoM}\right)\mathrm{I}_3-I_\mathrm{CoM}.\label{I_Sigma_IsTranslationInvariant}
\end{align}
This reveals a result similar to that of (\cite{wensing_linear_2018}, Proposition), specifically
\begin{align}
    \Sigma_\mathrm{CoM}\succ 0\quad\overset{\eqref{I_Sigma_IsTranslationInvariant}}{\Longleftrightarrow}\quad&\frac{1}{2}\mathrm{tr}(I_\mathrm{CoM})\mathrm{I}_3-I_\mathrm{CoM}\succ 0\allowdisplaybreaks\\
    \quad\Longleftrightarrow\quad&2\lambda_{\max}(I_\mathrm{CoM})<\mathrm{tr}(I_\mathrm{CoM})\allowdisplaybreaks\\
    \quad\Longleftrightarrow\quad&{\lambda_{\max}(I_\mathrm{CoM})<\lambda_{\min}(I_\mathrm{CoM})+\lambda_{\text{mid}}(I_\mathrm{CoM})}\allowdisplaybreaks\\
    \quad\Longleftrightarrow\quad&\underbrace{\lambda_{i}(I_\mathrm{CoM})<\lambda_{j}(I_\mathrm{CoM})+\lambda_{k}(I_\mathrm{CoM}) \text{ for all pairwise different } i, j, k\in\{1,2,3\}}_{\substack{\text{i.e.\ }I_\mathrm{CoM}\text{ fulfills the \textit{strict triangle inequalities} (see \cite{traversaro_identification_2016} for an introduction to their non-strict version)}}},\label{SigmaCoMPosDef_IFF_StricTriangleCoM}
\end{align}
wherein $\lambda_{\text{mid}}(I_\mathrm{CoM})=\mathrm{tr}(I_\mathrm{CoM})-\lambda_{\min}(I_\mathrm{CoM})-\lambda_{\max}(I_\mathrm{CoM})$ denotes the sandwiched eigenvalue of $I_\mathrm{CoM}$.

Necessity in the equivalence of physical consistency of $\Phi$ and positive definiteness of $f(\Phi)$ can be shown as in the proof of (\cite{wensing_linear_2018}, Theorem~$3$). For self-containment, we subsequently repeat their argument while utilizing our alternative~\eqref{SigmaCoMPosDef_IFF_StricTriangleCoM} to the required (\cite{wensing_linear_2018}, Proposition): Given $f(\Phi)\succ 0$, due to~\eqref{fPosDef_IFF_PosMassPosDefSigmaCoM}, we have $m>0$ and $\Sigma_\mathrm{CoM}\succ 0$ such that the sufficient direction in~\eqref{SigmaCoMPosDef_IFF_StricTriangleCoM} then guarantees $I_\mathrm{CoM}$ to fulfill the strict triangle inequalities. Thereby, physical consistency of $\Phi$, i.e.\ $m>0$ and realizability through some nonnegative mass density $\rho$ such that $\Phi$ satisfies~\eqref{MassDistributionToInertialParam}, follows by means of the construction of $\rho$ that is proposed in (\cite{traversaro_identification_2016}, Eq.\ $(26)$) since this $\rho$ satisfies~\eqref{MassDistributionToInertialParam} by design and maps exclusively into the nonnegative numbers whenever $m>0$ holds and $I_\mathrm{CoM}$ satisfies the non-strict triangle inequalities (“$\leq$" instead of “$<$" in the statement on right-hand-side of~\eqref{SigmaCoMPosDef_IFF_StricTriangleCoM}). As the already established fulfillment of the strict triangle inequalities of $I_\mathrm{CoM}$ trivially implies validity of their non-strict version, we know that $\rho$ as per (\cite{traversaro_identification_2016}, Eq.\ $(26)$) is nonnegative which, together with $m>0$, constitutes the desired physical consistency of $\Phi$.

The argument provided by \cite{wensing_linear_2018} for sufficiency in the equivalence of physical consistency of $\Phi$ and positive definiteness of $f(\Phi)$ can be summarized as follows: Given a physically consistent $\Phi$, one has existence of a nonnegative mass density $\rho$ that satisfies~\eqref{MassDistributionToInertialParam} with $m>0$. Therefore, as per (\cite{traversaro_identification_2016}, Eq.~(16)), one obtains that $I_\mathrm{CoM}$ fulfills the non-strict triangle inequalities via the nonnegativity of $\rho$. By means of the necessary directions in~\eqref{fPosDef_IFF_PosMassPosDefSigmaCoM} and in our alternative~\eqref{SigmaCoMPosDef_IFF_StricTriangleCoM} to (\cite{wensing_linear_2018}, Proposition),\footnote{To ensure that we fully capture the argument of~\cite{wensing_linear_2018} for sufficiency in the equivalence of physical consistency of $\Phi$ and positive definiteness of $f(\Phi)$, we subsequently explain the analogy of the required necessary directions in their (\cite{wensing_linear_2018}, Proposition) and in our~\eqref{SigmaCoMPosDef_IFF_StricTriangleCoM} since, for the sake of self containment, our discussion is facilitated by the latter. The former applied with respect to the pair $\Sigma_\mathrm{CoM}$, $I_\mathrm{CoM}$, that as per~\eqref{I_Sigma_IsTranslationInvariant} adheres to the requested structure $\Sigma_\mathrm{CoM}=\tfrac{1}{2}\mathrm{tr}(I_\mathrm{CoM})\mathrm{I}_3-I_\mathrm{CoM}$, does impose the same necessity as~\eqref{SigmaCoMPosDef_IFF_StricTriangleCoM} for the desired $\Sigma_\mathrm{CoM}\succ 0$, specifically that $I_\mathrm{CoM}$ needs to fulfill the strict triangle inequalities. However, the necessary direction in (\cite{wensing_linear_2018}, Proposition) additionally asks for positive definiteness of $I_\mathrm{CoM}$ in order to imply the desired $\Sigma_\mathrm{CoM}\succ 0$, which hints at a redundancy in its requirements on $I_\mathrm{CoM}$ since~\eqref{SigmaCoMPosDef_IFF_StricTriangleCoM} only asks for fulfillment of the strict triangle inequalities. This redundancy resolves as follows: Fulfillment of the strict triangle inequalities of $I_\mathrm{CoM}$ due to~\eqref{SigmaCoMPosDef_IFF_StricTriangleCoM} would imply $\lambda_k(\Sigma_\mathrm{CoM})>0$ for all $k\in\{1,2,3\}$, leading to $\mathrm{tr}(\Sigma_\mathrm{CoM})>\lambda_{\max}(\Sigma_\mathrm{CoM})$ such that $X=\mathrm{tr}(\Sigma_\mathrm{CoM})\mathrm{I}_3-\Sigma_\mathrm{CoM}\succ 0$. Note that taking the traces of both sides of~\eqref{I_Sigma_IsTranslationInvariant} yields $\mathrm{tr}(\Sigma_\mathrm{CoM})=\tfrac{1}{2}\mathrm{tr}(I_\mathrm{CoM})$ and thereby reveals $X=\tfrac{1}{2}\mathrm{tr}(I_\mathrm{CoM})-\Sigma_\mathrm{CoM}=I_\mathrm{CoM}$, cf.\ \eqref{I_Sigma_IsTranslationInvariant}. Therefore, the additional requirement $I_\mathrm{CoM}\succ 0$ for the necessary direction in (\cite{wensing_linear_2018}, Proposition) to hold can be dropped since its validity is implied whenever $I_\mathrm{CoM}$ satisfies the strict triangle inequalities, i.e.\ whenever the desired $\Sigma_\mathrm{CoM}\succ 0$ is guaranteed by means of the necessary direction in~\eqref{SigmaCoMPosDef_IFF_StricTriangleCoM}. Hence, the conditions for the necessary directions in~\eqref{SigmaCoMPosDef_IFF_StricTriangleCoM} and (\cite{wensing_linear_2018}, Proposition) are equally demanding and thus can be used interchangeably.} the desired $f(\Phi)\succ 0$ then would follow from $\Sigma_\mathrm{CoM}\succ 0$ and $m>0$ if one were able to rule out equality in all three triangle inequalities of $I_\mathrm{CoM}$,~i.e.\ by ensuring their strictness, which, even though illustrated with some example point-mass distributions, is not proven in \cite{wensing_linear_2018} for arbitrary rigid bodies with physically consistent inertial parameters; thus leaving a gap in the reasoning for why $f(\Phi)\succ 0$ holds under physical consistency of $\Phi$.

Following the intuition built by the illustrative examples given in \cite{wensing_linear_2018} that strictness in the triangle inequalities of $I_\mathrm{CoM}$ encodes the $3$-dimensionality of a rigid body, we choose to fill the aforementioned gap in the theory by exploiting properties of the shape $\{x\in\mathbb{R}^3:\rho(x)\neq 0\}$ of a rigid body, specifically that its measurability---under the presupposition that the integrals in~\eqref{MassDistributionToInertialParam} are evaluated with respect to the Lebesgue measure for volumina in $\mathbb{R}^3$---is closely related to physical consistency of the inertial parameters $\Phi$. In order to exploit this connection, we  restrict the integrals in~\eqref{MassDistributionToInertialParam} to Lebesgue integration, thereby excluding any mass densities $\rho$ defined by distributions such as point masses, line masses or planar masses from the scope of our analysis, and thus arrive at the refined definition of physical consistency~in: 

\begin{definition}\label{definition:PhysicalConsistency}
    A vector $\Phi\in\mathbb{R}^{10}$ with entries as in $(\ref{eq:Phi})$ is called \textit{physically consistent} if there exists a nonnegative mass density, specifically a function $\rho:\mathbb{R}^3\to[0,\infty)$, such that the equations in~\eqref{MassDistributionToInertialParam} hold with $m>0$ when the integrals therein are evaluated with respect to the Lebesgue measure for volumina in $\mathbb{R}^3$.
\end{definition}

The physical consistency of inertial parameters then is exactly characterized by the next statement.

\begin{theorem}\label{theorem:PhysicalConsistency} A vector $\Phi\in\mathbb{R}^{10}$ is physically consistent (in accordance with the refined  Definition~\ref{definition:PhysicalConsistency}) if and only if $f(\Phi)\succ 0$. Furthermore, it is noteworthy that the inverse transformation from the upper left block $\Sigma$ of $f(\Phi)$ to the components of the inertia matrix $I$ that build the vector $\Phi$ reads as $I=\mathrm{tr}(\Sigma)\mathrm{I}_3-\Sigma$.
\end{theorem}

\begin{proof}
    First, the inverse transformation follows analogously to (\cite{wensing_linear_2018}, Eqs.~(16),~(17)) by equating the traces of both sides of $\Sigma=\tfrac{1}{2}\mathrm{tr}(I)\mathrm{I}_3-I$. Second, necessity in the equivalence claimed in Theorem~\ref{theorem:PhysicalConsistency} is already shown with the discussion above since the mass density (\cite{traversaro_identification_2016}, Eq.\ $(26)$) is a nonnegative realization of any physically consistent $\Phi$ when using Lebesgue integration in~\eqref{MassDistributionToInertialParam}. The proof of sufficiency in the equivalence claimed in Theorem~\ref{theorem:PhysicalConsistency} is moved to the Appendix for better readability. There, we prove the sufficient direction by removing the possibility of semi-definiteness of $f(\Phi)$ associated with equality in the triangle inequalities of $I_\mathrm{CoM}$ via application of the strict version of the Schwarz inequality. Crucially, we exploit the restriction in Definition~\ref{definition:PhysicalConsistency} that mass densities defining physically consistent rigid bodies must be nonnegative \textit{functions} with $m>0$, which guarantees positive measure of $\big\{x\in\mathbb{R}^3:\rho(x)\neq 0\big\}$ and thus excludes flat bodies contained in a hyperplane that lack the $3$-dimensionality required for strictness of the Schwarz inequality which, under $m>0$, can be shown to be equivalent to the desired $f(\Phi)\succ 0$.
\end{proof}

\begin{remark}
The intrinsic relationship~\eqref{MassDistributionToInertialParam} is the reason why, when the rigidity of a body is relaxed, time dependency of $\Phi$ necessitates time-dependent change of the mass density $\rho$. As time dependency of $\rho$ is also triggered by movements of mass-carrying particles relative to their body fixed frame, the effects of these movements are included in the model that we derive in Section~\ref{section:Modeling}.
\end{remark}

\section{Generalized robotics equation}\label{section:Modeling}
We consider a robotic manipulator that consists of an $n$-degree-of-freedom open kinematic chain of $N$ bodies. In order to describe its dynamical behavior, first, we introduce the model and then the kinematics. We discuss the restrictions imposed by the modeling that are essential for obtaining an ODE generalization of~\eqref{SimpleGeneralizationOfTheRoboticEquation} from \cite{pagilla_adaptive_2000}. This generalization then captures all relevant effects of time-dependent inertial parameters as well as the effects of the causative time dependencies of the mass densities originating from internal movements of mass-carrying particles. However, by means of the standing Assumption~\ref{Assumption:IndependentJacobian}, we exclude other effects that would result from allowing time dependency of kinematic parameters. Subsequently, we derive the kinetic energies $\mathcal{T}_l$ as well as the potential energies $\mathcal{U}_l$ of the individual bodies $l\in\{1,\ldots,N\}$. Finally, these energies are reformulated with respect to the specific movement of the robotic manipulator governed by its forward kinematic map and the corresponding stacked Jacobian so that we then obtain the generalized robotics equation describing the dynamical behavior by means of the Lagrange formalism.

\paragraph{Model} The pose of the $l$th body-fixed frame is given through the position $z_{l}(t)\in\mathbb{R}^3$ of its origin in the $0$th frame fixed to a Newtonian system and by the rotation angles $\phi_l(t)\in\mathbb{R}^3$ that parameterize its orientation with the rotation matrix $\mathcal{R}(\phi_l)\in\mathbb{R}^{3\times 3}$ such that any point $x_l\in\mathbb{R}^3$ in the $l$th body-fixed frame reads as $\mathcal{R}(\phi_l)x_l+z_l$ in the $0$th Newtonian frame. All frames under consideration are equipped with orthonormal bases, which renders the rotation matrix orthogonal, i.e.\ $\mathcal{R}(\phi_l)^\top\mathcal{R}(\phi_l)=\mathrm{I}_3$. Thereby, $\mathcal{R}(\phi_l)\in\mathrm{SO}(3)=\{M\in\mathbb{R}^{3\times 3}: \mathrm{det}(M)=1,\ M^\top=M^{-1}\}$ as a rotation inherently satisfies $\mathrm{det}(R(\phi_l))>0$. Further, the translational velocity of the origin of the $l$th body-fixed frame is denoted as $v_l(t)=\dot{z}_l(t)\in\mathbb{R}^3$ and $\omega_l(t)\in\mathbb{R}^3$ is its angular velocity with coordinates expressed in the $l$th body-fixed frame, i.e.~$\mathcal{S}(\omega_l)=\mathcal{R}(\phi_l)^\top \dot{\mathcal{R}}(\phi_l)$. Let $\rho_l(x_l,t)\in\mathbb{R}$ the respective density of mass-carrying particles and $\sigma_l(x_l,t)\in[0,1]$ the portion of these particles moving with velocity $\mathscr{v}_l(x_l,t)\in\mathbb{R}^3$  relative to the $l$th body-fixed frame at the position $x_l\in\mathbb{R}^3$ in this frame. The remaining portion $1-\sigma_l(x_l,t)$ of mass-carrying particles at this position is immobile with respect to the frame that is fixed to the $l$th body. In accordance with~\eqref{eq:Phi} and~\eqref{MassDistributionToInertialParam}, we define the consequently time-dependent inertial parameters $\Phi_l(t)\in\mathbb{R}^{10}$ of the $l$th body as
\begin{align}
    \Phi_l&=\left(\begin{matrix}m_l  & h_{1,l}  & h_{2,l}  & h_{3,l} & I_{11,l}  & I_{22,l}  & I_{33,l}  & I_{12,l}  & I_{23,l}  & I_{13,l}\end{matrix}\right)^\top,
\end{align}
where
\begin{align}
\label{Defml}
        m_l(\cdot)&=\int_{\mathbb{R}^3} \rho_l(x_l,\cdot)\mathrm{d}x_l\in\mathbb{R},\allowdisplaybreaks\\
        h_l(\cdot)&=\begin{pmatrix}h_{1,l}(\cdot)&h_{2,l}(\cdot)&h_{3,l}(\cdot)\end{pmatrix}^\top=\int_{\mathbb{R}^3} \rho_l(x_l,\cdot)x_l\mathrm{d}x_l\in\mathbb{R}^3,\allowdisplaybreaks\\
        I_l(\cdot)&=\left(\begin{smallmatrix}
        I_{11,l}(\cdot)&I_{12,l}(\cdot)&I_{13,l}(\cdot)\\
        I_{12,l}(\cdot)&I_{22,l}(\cdot)&I_{23,l}(\cdot)\\
        I_{13,l}(\cdot)&I_{23,l}(\cdot)&I_{33,l}(\cdot)
    \end{smallmatrix}\right)=\int_{\mathbb{R}^3} \rho_l(x_l,\cdot)\mathcal{S}(x_l)^\top \mathcal{S}(x_l)\mathrm{d}x_l\in\mathrm{sym}(3).\label{DefIl}
\end{align}
The vector $\Theta(t)=\left(\Theta_h(t)\right)_{h\in\{1,\ldots,10N\}}\in\mathbb{R}^{10N}$ with
\begin{align}
\Theta=\begin{pmatrix}\Phi_1^\top&\ldots&\Phi_N^\top\end{pmatrix}^\top    
\end{align}
collects the complete set of inertial parameters describing the robotic manipulator.
\paragraph{Discussion of restrictions imposed by the model}
    The model from above includes time-dependent change of the mass densities and thus time dependency of the inertial parameters by means of
    \begin{itemize}
        \item internal mass-density redistribution represented by $\mathscr{v}_l(x_l,t)\neq 0$ and $\sigma_l(x_l,t)\in(0,1]$, or
        \item the addition of mass via $\dot{\rho}_l(x_l,t)\neq 0$
    \end{itemize}
     for some $x_l\in\mathbb{R}^3$, $t\geq 0$. However, in order to ensure that the dynamical behavior results as ODE, not as PDE, the distributed quantities $\mathscr{v}_l$, $\sigma_l$ and $\dot{\rho}_l$ that determine the movements of mass-carrying particles relative to their body-fixed frames are modeled as being independent of $z_l,\phi_l$ or their time derivatives. This excludes effects like elasticity of the bodies from the model since the time dependency is seen as an external process that might affect the dynamical behavior but is not driven by the movement of the robotic manipulator.

\paragraph{Kinematics} The kinematics describe how the $N$ bodies interact when their movements are subject to restrictions imposed by the joints of the robotic manipulator. To capture this effect, we introduce the generalized coordinates $q(t)=\left(q_k(t)\right)_{k\in\{1,\ldots,n\}}\in\mathbb{R}^n$. Since the kinematic chain is \textit{open}, $q(t)$ is a  \textit{minimal set of pairwise independent variables} such that a forward kinematic map $F(q,\Theta_\mathrm{kin})\in\mathbb{R}^{6N}$ with
\begin{align}
\label{Layout}
    \begin{pmatrix}
         z_1^\top&\phi_1^\top&\ldots&z_N^\top&\phi_N^\top
    \end{pmatrix}^\top&=F(q,\Theta_\mathrm{kin})
\end{align}
exists for all poses $z_l, \phi_l\in\mathbb{R}^{3}$, $l\in\{1,\ldots,N\}$ that are possible during operation of the robotic manipulator. Therein, $\Theta_\mathrm{kin}\in\mathbb{R}^{n_\mathrm{kin}}$ are kinematic parameters (lengths and angles describing the poses of the joints in the respective body-fixed frames) defining the kinematic behavior of the robotic manipulator. In general, due to effects like thermal expansion, the kinematic parameters depend on the mass densities. That is, there exists a function $D(\rho_1,\ldots,\rho_N)\in\mathbb{R}^{n_\mathrm{kin}}$ with
\begin{align}
\Theta_\mathrm{kin}=&D(\rho_1,\ldots,\rho_N)
\end{align}
such that, in general, the kinematic parameters inherit time dependency from the mass densities. However---since our main goal is to highlight the effect of time-dependent inertial parameters as well as the causative time dependency of the mass densities---we restrict the analysis in this contribution to robotic manipulators that satisfy the following standing assumption:
\begin{assumption}\label{Assumption:IndependentJacobian}
    The function $D(\rho_1,\ldots,\rho_N)$ is constant, i.e.\ the kinematic parameters remain constant despite possible time-dependent change of the mass densities $\rho_1,\ldots,\rho_N$ and despite the resulting time dependency of the inertial parameters in $\Theta$. Accordingly, we have
    \begin{align}
    \label{Kinematics}
    {\begin{pmatrix}
        v_1^\top&\dot{\phi}_1^\top&\ldots&v_N^\top&\dot{\phi}_N^\top
    \end{pmatrix}^{\top}}=J(q)\dot{q}
    \end{align}
    with the Jacobian\footnote{Please note that here in this work, $J$ is the \textit{stacked} Jacobian of the forward kinematic map $F$ that relates the generalized coordinates to the movement of \textit{all} bodies in the robotic manipulator; a concept that is distinct from and should not be confused with that of an end-effector Jacobian $J_{l_\mathrm{ee}\!}(q)=\tfrac{\partial}{\partial q}(z_{l_\mathrm{ee}\!}\!^\top,\phi_{l_\mathrm{ee}\!}\!^\top)\!^\top$, where $l_{\mathrm{ee}\!}\in\{1,\ldots,N\}$ is an index that refers to a single end-effector.}
    \begin{align}
    J(q)=\tfrac{\partial}{\partial q}F(q,\Theta_\mathrm{kin})
    \end{align}
    of the forward kinematic map $F$, wherein we have dropped the dependency on the kinematic parameters for better readability.
\end{assumption}

\paragraph{Discussion of the standing Assumption~\ref{Assumption:IndependentJacobian}.} In order to understand the prerequisites for the standing Assumption~\ref{Assumption:IndependentJacobian}, consider the following: Usually, the origins of the body-fixed frames are chosen at the joints between the bodies. Then, Assumption~\ref{Assumption:IndependentJacobian} is fulfilled if the   forward kinematic map $F$ and therewith the poses of the joints relative to each other are independent of the mass densities of the bodies that are linking them. That is, time dependency of the inertial parameters while Assumption~\ref{Assumption:IndependentJacobian} holds and with the body-fixed frames located at the joints, can be understood as the effect of time-dependent change of the mass densities of the links between the joints that neither affects the positions nor the orientations of the joints relative to each other. As a remarkable consequence of this, provided that the body-fixed frames are located at the joints, arbitrary change with time of the mass density of an end-effector does not interfere with the validity of the standing Assumption~\ref{Assumption:IndependentJacobian}, since end-effectors do not serve as links between joints. 

\paragraph{Calculation of the energies.} Subsequently, we aim at calculating the kinetic energies of the individual bodies. To that aim, we need the norms of the velocities with respect to the $0$th Newtonian frame of the particles that build these bodies. For a particle that is stationary at $x_l$ in its body fixed frame, the norm of its velocity in the Newtonian system is
\begin{align}
    \left\|\frac{\mathrm{d}}{\mathrm{d}t}(z_l+\mathcal{R}(\phi_l)x_l)\right\|=&\left\|v_l+\dot{\mathcal{R}}(\phi_l)x_l\right\|=\left\|\mathcal{R}(\phi_l)\left(\mathcal{R}(\phi_l)^\top v_l+\mathcal{R}(\phi_l)^\top \dot{\mathcal{R}}(\phi_l)x_l\right)\right\|\allowdisplaybreaks\\
    \label{VelocityRigidMass}
    =&\left\|\mathcal{R}(\phi_l)^\top v_l+\mathcal{S}(\omega_l)x_l\right\|.
\end{align}
Note that the rotation matrix allows calculating the velocities of mass-carrying particles relative to their body-fixed frames in the $0$th Newtonian frame as $\mathcal{R}(\phi_l)\mathscr{v}_l$. Accordingly, for a particle located at $x_l$ in its body-fixed frame that moves with velocity $\mathscr{v}_l(x_l,t)$ relative to this frame, the norm of its velocity in the Newtonian system is
\begin{align}
    \left\|\mathcal{R}(\phi_l)\mathscr{v}_l(x_l,t)+\frac{\mathrm{d}}{\mathrm{d}t}(z_l+\mathcal{R}(\phi_l)x_l)\right\|=&\left\|\mathcal{R}(\phi_l)\left(\mathscr{v}_l(x_l,t)+\mathcal{R}(\phi_l)^\top v_l+\mathcal{R}(\phi_l)^\top \dot{\mathcal{R}}(\phi_l)x_l\right)\right\|\allowdisplaybreaks\\\label{VelocitySquishyMass}
    =&\left\|\mathscr{v}_l(x_l,t)+\mathcal{R}(\phi_l)^\top v_l+\mathcal{S}(\omega_l)x_l\right\|.
\end{align}
The velocities~\eqref{VelocityRigidMass},~\eqref{VelocitySquishyMass} relative to the Newtonian system of particles that are fixed and of those that move inside their body-fixed frame, respectively, both contribute to the kinetic energy of the $l$th body, which thereby amounts to
\begin{align}
    \mathcal{T}_l&=\frac{1}{2}\int_{\mathbb{R}^3}\left((1-\sigma_l(x_l,t))\rho_l(x_l,t)\left\|\mathcal{R}(\phi_l)^\top v_l+\omega_l\times x_l\right\|^2\right.\\
    &\quad \left.+\sigma_l(x_l,t)\rho_l(x_l,t)\left\|\mathscr{v}_l(x_l,t)+\mathcal{R}(\phi_l)^\top v_l+\omega_l\times x_l\right\|^2\right)\mathrm{d}x_l\nonumber\allowdisplaybreaks\\
    &=\frac{1}{2}\int_{\mathbb{R}^3}\left((1-\sigma_l(x_l,t))\rho_l(x_l,t)\left\|\mathcal{R}(\phi_l)^\top v_l-\mathcal{S}(x_l)\omega_l\right\|^2\right.\\ &\quad \left.+\sigma_l(x_l,t)\rho_l(x_l,t)\left\|\mathscr{v}_l(x_l,t)+\mathcal{R}(\phi_l)^\top v_l-\mathcal{S}(x_l)\omega_l\right\|^2\right)\mathrm{d}x_l\nonumber
    \allowdisplaybreaks\\
    &=\frac{1}{2}\int_{\mathbb{R}^3} \left(\rho_l(x_l,t) \left(\left\|\mathcal{R}(\phi_l)^\top v_l\right\|^2-2v_l^\top\mathcal{R}(\phi_l) \mathcal{S}(x_l)\omega_l+ \left\|\mathcal{S}(x_l)\omega_l\right\|^2\right)\right.\\
    &\quad \left.+ \sigma_l(x_l,t)\rho_l(x_l,t)\left(\|\mathscr{v}_l(x_l,t)\|^2+2\mathscr{v}_l(x_l,t)^\top (\mathcal{R}(\phi_l)^\top v_l-\mathcal{S}(x_l)\omega_l)\right)\right)\mathrm{d}x_l\nonumber\allowdisplaybreaks\\
    \label{RespectiveKineticEnergy_Intermediate}
    &=\frac{1}{2}\int_{\mathbb{R}^3}\rho_l(x_l,t)\begin{pmatrix}
        \mathcal{R}(\phi_l)^\top v_l\\\omega_l
    \end{pmatrix}^{\top} \begin{pmatrix}
        \mathrm{I}_3&-\mathcal{S}(x_l)\\-\mathcal{S}(x_l)^\top&\mathcal{S}(x_l)^\top \mathcal{S}(x_l)
    \end{pmatrix}\begin{pmatrix}
        \mathcal{R}(\phi_l)^\top v_l\\\omega_l
    \end{pmatrix}\mathrm{d}x_l\\
    &\quad +\frac{1}{2}\int_{\mathbb{R}^3}  \sigma_l(x_l,t)\rho_l(x_l,t)\|\mathscr{v}_l(x_l,t)\|^2   \mathrm{d}x_l+\int_{\mathbb{R}^3}\sigma_l(x_l,t)\rho_l(x_l,t)\begin{pmatrix}
        \mathcal{R}(\phi_l)^\top v_l\\\omega_l
    \end{pmatrix}^{  \top} \begin{pmatrix}
        \mathrm{I}_3\\\mathcal{S}(x_l)
    \end{pmatrix}\mathscr{v}_l(x_l,t)\mathrm{d}x_l\nonumber.
\end{align}
In view of~\eqref{Defml}-\eqref{DefIl}, this leads to
\begin{align}
    \label{IntermediateKineticEnergy}
    \mathcal{T}_l&=\frac{1}{2}\begin{pmatrix}
        \mathcal{R}(\phi_l)^\top v_l\\\omega_l
    \end{pmatrix}^{\top} Z(\Phi_l)\begin{pmatrix}
        \mathcal{R}(\phi_l)^\top v_l\\\omega_l
    \end{pmatrix}+\frac{1}{2}\nu_l(t) +\begin{pmatrix}
        \mathcal{R}(\phi_l)^\top v_l\\\omega_l
    \end{pmatrix}^{\top} \psi_l(t)
\end{align}
for $l\in\{1,\ldots,N\}$ with
\begin{align}
\label{DefMatrixZ}
    Z(\Phi_l)&=\begin{pmatrix}
        m_l\mathrm{I}_3&-\mathcal{S}(h_l)\\
        -\mathcal{S}(h_l)^\top &I_l
    \end{pmatrix}\in\mathrm{sym}(6),\allowdisplaybreaks\\
    \nu_l(t)&=\int_{\mathbb{R}^3}\sigma_l(x_l,t)\rho_l(x_l,t)\left\|\mathscr{v}_l(x_l,t)\right\|^2\mathrm{d}x_l\in\mathbb{R},\allowdisplaybreaks\\
    \psi_l(t)&=\int_{\mathbb{R}^3}\sigma_l(x_l,t)\rho_l(x_l,t)\begin{pmatrix}
       \mathrm{I}_3\\\mathcal{S}(x_l)
    \end{pmatrix}\mathscr{v}_l(x_l,t)\mathrm{d}x_l\in\mathbb{R}^{6}.
\end{align}
Further, due to $\mathcal{S}(\omega_l)=\mathcal{R}(\phi_l)^\top \dot{\mathcal{R}}(\phi_l)$, there exists a map $E(\phi_l)\in\mathbb{R}^{3\times 3}$ with $ \omega_l=E(\phi_l)\dot{\phi}_l$.
Therefore, the kinetic energy of the $l$th body in~\eqref{IntermediateKineticEnergy} rewrites as
\begin{align}
    \label{RespectiveKineticEnergy}
    \mathcal{T}_l&=\frac{1}{2}\begin{pmatrix}
        v_l\\\dot{\phi}_l
    \end{pmatrix}^{\top}Q(\phi_l)^\top Z(\Phi_l)Q(\phi_l)\begin{pmatrix}
        v_l\\\dot{\phi}_l
    \end{pmatrix}+\frac{1}{2}\nu_l(t) +\begin{pmatrix}
        v_l\\\dot{\phi}_l
    \end{pmatrix}^{\top}Q(\phi_l)^\top \psi_l(t)
\end{align}
with the mapping
\begin{align}
    Q(\phi_l)=&\mathrm{diag}\left(\mathcal{R}(\phi_l)^\top, E(\phi_l)\right)\in\mathbb{R}^{6\times 6}
\end{align}
that describes the specific method for parameterizing the orientations of the body-fixed frames with the rotation angles.  Next, by presupposing that the robotic manipulator operates in a constant gravitational field $g\in\mathbb{R}^3$ with respect to the $0$th Newtonian frame, we receive the potential energies of the individual bodies $l\in\{1,\ldots,N\}$ as
\begin{align}
    \mathcal{U}_l&=m_lg^\top(z_l+\mathcal{R}(\phi_l)c_l)\allowdisplaybreaks\\
    \label{RespectivePotenitalEnergy}
    &= g^\top(m_l z_l+ \mathcal{R}(\phi_l)h_l),
\end{align}
where $c_l(t)=m_l(t)^{-1}h_l(t)\in\mathbb{R}^3$ is the location of the CoM of the $l$th rigid body in its body-fixed frame.

We proceed by reformulating these kinetic and potential energies with respect to the specific movement of the bodies possible during operation of the robotic manipulator by using its kinematics. By means of~\eqref{RespectiveKineticEnergy} and~\eqref{Kinematics}, the overall kinetic energy of the open kinematic chain can be expressed in dependence of the generalized coordinates, i.e.\ with respect to the movement of the robotic manipulator, as
\begin{align}
        \mathcal{T}(q,\dot{q},\Theta,t)&=\sum_{l=1}^N \mathcal{T}_l
        \label{KineticEnergy}=\frac{1}{2}\dot{q}^\top M(q,\Theta)\dot{q}+\frac{1}{2}\nu(t)+\dot{q}^\top J(q)^\top \mathcal{Q}(q)^\top \Psi(t),
\end{align}
where the mass matrix results as
\begin{align}
    \label{MassMatrix}
M(q,\Theta)&=J(q)^\top \mathcal{Q}(q)^\top \mathcal{Z}(\Theta)\mathcal{Q}(q)J(q)
\end{align}
with
\begin{align}
\label{DefBlockmatrixZ}
    \mathcal{Z}(\Theta)&=\mathrm{diag} \left(
    Z(\Phi_1), \ldots, Z(\Phi_N)
\right)\in\mathrm{sym}(6N)\allowdisplaybreaks,\\
\label{DefBlockmatrixQ}
\mathcal{Q}(q)&=\mathrm{diag} \left(
    Q\left(F_{\mathrm{\phi}_1}(q)\right), \ldots, Q\left(F_{\mathrm{\phi}_N}(q)\right)
\right)\in\mathbb{R}^{6N\times 6N}.
\end{align}
Therein, $F_{\mathrm{z}_l}(q),F_{\mathrm{\phi}_l}(q)\in\mathbb{R}^3$, $l\in\{1,\ldots,N\}$ represent a decomposition of the forward kinematic map such that $F(q,\Theta_\mathrm{kin})=(\begin{smallmatrix}
    F_{\mathrm{z}_1}(q)^\top&F_{\mathrm{\phi}_1}(q)^\top&\ldots&F_{\mathrm{z}_N}(q)^\top&F_{\mathrm{\phi}_N}(q)^\top
\end{smallmatrix}) ^\top$. The remaining terms in~\eqref{KineticEnergy} are
\begin{align}
\nu(t)&=\sum_{l=1}^N \nu_l(t)\in\mathbb{R},\allowdisplaybreaks\\
\Psi(t)&=\begin{pmatrix}
 \psi_1(t)^\top&\ldots&\psi_N(t)^\top
\end{pmatrix}^\top\in\mathbb{R}^{6N}.
\end{align}
The effects of the latter two contributions in the overall kinetic energy~\eqref{KineticEnergy} associated with $\nu$ and $\Psi$ are absent in the version of the robotics equation from~\cite{pagilla_adaptive_2000} and thus, they represent the structural novelty of our model, which is the inclusion of internal mass-density redistribution. It should be noted, however, that only the lumped velocity $\Psi$ and its time derivative $\dot{\Psi}$ will appear in the generalized robotics equation below since the energy contribution associated with $\nu$ is independent of the generalized coordinates $q$, independent of their time derivatives $\dot{q}$ due to our modeling choice to understand the internal movement of mass-carrying particles as an external process in order to keep the dynamical behavior as an ODE. Further, note that $M(q,\Theta)=\big(M_{k,j}(q,\Theta)\big)_{(k,j)\in\{1,\ldots,n\}^2}\in\mathrm{sym}(n)$ is linearly dependent on the inertial parameters in $\Theta$ because they appear linearly in the block matrix $\mathcal{Z}(\Theta)$.

In similar fashion to the derivation of the overall kinetic energy, due to~\eqref{RespectivePotenitalEnergy} and~\eqref{Layout}, the overall potential energy of the open kinematic chain in dependence of the generalized coordinates, i.e.\ with respect to the movement of the robotic manipulator, evaluates to
\begin{align}
    \mathcal{U}(q,\Theta)&=\sum_{l=1}^N\mathcal{U}_l=g^\top\Big(\sum_{l=1}^N m_l F_{\mathrm{z}_l}(q)+ \mathcal{R} \left(F_{\mathrm{\phi}_l}(q)\right)h_l\Big),
    \label{PotentialEnergy}
\end{align}
which is also linearly dependent on the inertial parameters in $\Theta$, as is evident by its structure.

\paragraph{Lagrange formalism} In~\eqref{KineticEnergy} and~\eqref{PotentialEnergy}, we have gathered enough information to write the Lagrangian
\begin{align}
\label{Lagrangian}
    L(q,\dot{q},\Theta,t)&=\mathcal{T}(q,\dot{q},\Theta,t)-\mathcal{U}(q,\Theta)
\end{align}
which, by virtue of the Lagrange formalism, leads to the nonlinear dynamical behavior
\begin{align}
\label{LagrangeFormalism}
    \frac{\mathrm{d}}{\mathrm{d}t}\left(\frac{\partial L(q,\dot{q},\Theta,t)}{\partial \dot{q}_k}\right)-\frac{\partial L(q,\dot{q},\Theta,t)}{\partial q_k}=\tau_k+w_k,
\end{align}
$k\in\{1,\ldots,n\}$ of the robotic manipulator, where $\tau(t)=\big(\tau_k(t)\big)_{k\in\{1,\ldots,n\}}\in\mathbb{R}^n$ is the external torque/force that gets applied at the joints together with some disturbance $w(t)=\big(w_k(t)\big)_{k\in\{1,\ldots,n\}}\in\mathbb{R}^n$. Now, as is shown in the Appendix, the dynamical behavior~\eqref{LagrangeFormalism} rewrites as \textbf{generalized robotics equation} 
\begin{align}
\label{RoboticEquation}
    &M(q,\Theta)\ddot{q} + \left(C(q,\dot{q},\Theta)+M(q,\dot{\Theta})+H(q,\Psi)\right)\dot{q}+G(q,\Theta)=\tau+w-J(q)^\top \mathcal{Q}(q)^\top\dot{\Psi}.
\end{align}
The matrices in~\eqref{RoboticEquation} are
\begin{align}
\label{Cristoffel}
    \Gamma_{i,j,k}(q,\Theta)&=\frac{1}{2}\left(\frac{\partial M_{k,j}(q,\Theta)}{\partial q_i}+\frac{\partial M_{k,i}(q,\Theta)}{\partial q_j}-\frac{\partial M_{i,j}(q,\Theta)}{\partial q_k}\right)\in\mathbb{R},\allowdisplaybreaks\\
\label{DefCoriolis}
C(q,\dot{q},\Theta)&=\left(\sum_{i=1}^n\Gamma_{i,j,k}(q,\Theta)\dot{q}_i\right)_{(k,j)\in\{1,\ldots,n\}^2}\in\mathbb{R}^{n\times n},\allowdisplaybreaks\\
\label{DefGravitation}
G(q,\Theta)&=\left(\frac{\partial U(q,\Theta)}{\partial q}\right)^{\top}\in\mathbb{R}^n,\allowdisplaybreaks\\
\label{DefH}
H(q,\Psi)&=\frac{\partial J(q)^\top \mathcal{Q}(q)^\top \Psi}{\partial q}-\left(\frac{\partial J(q)^\top \mathcal{Q}(q)^\top \Psi}{\partial q}\right)^{ \top}\in\mathbb{R}^{n\times n}
\end{align}
with the Christoffel symbols of the first kind $\Gamma_{i,j,k}(q,\Theta)$, $i,j,k\in\{1,\ldots,n\}$ building the Coriolis matrix $C(q,\dot{q},\Theta)$ and $G(q,\Theta)$ representing the influence of gravitation, as is already known from the classical robotics equation~\eqref{ClassicalRoboticEquation}. However, the generalized robotics equation~\eqref{RoboticEquation} also includes some effects beyond the scope of its classical counterpart:
\begin{itemize}
\item The effect of time dependency of the inertial parameters in $\Theta$ on the dynamical behavior is expressed in~\eqref{RoboticEquation} through the contribution that is made by $M(q,\dot{\Theta})\dot{q}$. This effect is already represented with $\mathcal{F}(q,\Theta,\dot{\Theta})\dot{q}$ in the version of the robotics equation~\eqref{SimpleGeneralizationOfTheRoboticEquation} from~\cite{pagilla_adaptive_2000}. However, the derivation of the generalized robotics equation in the Appendix reveals $\mathcal{F}(q,\Theta,\dot{\Theta})=M(q,\dot{\Theta})$. The description with $M(q,\dot{\Theta})$ offers more structural insight that we make explicit with Theorem~\ref{theorem:BoundedRateOfChange} in Section~\ref{section:OtherComponents}.
\item The distributed velocities $\mathscr{v}_l$, $l\in\{1,\ldots,N\}$ of mass-carrying particles relative to their body-fixed frames are lumped into $\Psi$ and their effect on the dynamical behavior appears in~\eqref{RoboticEquation} via $H(q,\Psi)\dot{q}$.
\item The force exerted by mass-carrying particles that experience acceleration $\dot{\Psi}$ relative to their body-fixed frames in a lumped sense is represented in~\eqref{RoboticEquation} by $J(q)^\top\mathcal{Q}(q)^\top\dot{\Psi}$ and acts like a disturbance.
\end{itemize}

\section{Structural properties}\label{section:StructuralProperties}

In this section, we exploit insight from the modeling in order to reveal structural properties of the generalized robotics equation~\eqref{RoboticEquation} and of some of its components. In principle, the following two properties inherent to robotics equations are already known, see e.g.\ \cite{pagilla_adaptive_2000}. However for the sake of completeness, we state them here in a form tailored to our generalization of the robotics equation. Their derivation is moved to the Appendix and strictly adheres to the model as formulated in Section~\ref{section:Modeling}. We have:
\begin{enumerate}
    \item\label{Prop1} The matrices $\dot{M}(q,\Theta)-2C(q,\dot{q},\Theta)-M(q,\dot{\Theta})$ and $H(q,\Psi)$ are skew-symmetric for all $q,\dot{q}\in \mathbb{R}^{n}$, $\Theta,\dot{\Theta}\in\mathbb{R}^{10N}$, $\Psi\in\mathbb{R}^{6N}$.
    \item\label{Prop2} The generalized robotics equation~\eqref{RoboticEquation} is linearly dependent on constant inertial parameters, i.e.\ there exist regressor functions $R_l(q,\dot{q},v,a)\in\mathbb{R}^{n\times 10}$,  $l\in\{1,\ldots,N \}$ with
    \begin{align}
    \label{RoboticRegressor}
        M (q,\Theta)a + C (q,\dot{q},\Theta)v + G (q,\Theta)
        = &\sum_{l=1}^N R_l (q,\dot{q},v,a)\Phi_l
    \end{align}
    for all $q,\dot{q},v,a\in\mathbb{R}^n$, $\Phi_1,\ldots,\Phi_N\in\mathbb{R}^{10}$. Furthermore, there also exist regressors $R_{\mathrm{M},l}(q,v)\in\mathbb{R}^{n\times 10}$, $l\in\{1,\ldots,N \}$ that allow representing the effect of time dependency of inertial parameters as
    \begin{align}
    \label{RoboticRegressorVel}
        M (q,\dot{\Theta})v
        = &\sum_{l=1}^N R_{\mathrm{M},l}(q,v)\dot{\Phi}_l,
    \end{align}
    which holds for all $q,v\in\mathbb{R}^n$, $\dot{\Phi}_1,\ldots,\dot{\Phi}_N\in\mathbb{R}^{10}$.
\end{enumerate}
The next two subsections are concerned with the derivation of conditions under which components of the generalized robotics equation~\eqref{RoboticEquation} in some sense are bounded. This is done with the intention to enable analysis for robustness of adaptive algorithms as laid out in the introduction to conclude (ultimate) boundedness of signals in the closed loop when unknown inertial parameters are at risk to depend on time in a manner unforeseen by the adaptation but these conditions are guaranteed to hold.

\subsection{Boundedness of the mass matrix}\label{section:BoundedMassMatrix}

The tracking capabilities of passivity-based adaptive control methods for robotic manipulators such as \cite{slotine_adaptive_1987}, its successor \cite{johansson_adaptive_1990} or the approach in \cite{pagilla_adaptive_2000} are guaranteed by means of Lyapunov functions that in part are formed as a quadratic form of the tracking error with the mass matrix. Therefore, the guarantees of stability of such control algorithms hinge on the existence of a positive lower uniform bound and a finite upper uniform bound of the mass matrix, thus rendering conditions for the existence of such finite, positive uniform bounds desirable. Before we highlight challenges arising explicitly from time dependency of the inertial parameters, we provide some insight into the existence of such bounds in the case of constant, physically consistent inertial parameters: First of all, positivity of a lower uniform bound of the mass matrix necessitates the property in

\begin{definition}\label{definition:NormalJacobian}
    $\mathcal{Q}(\cdot)J(\cdot)$ is called \textit{normal} if $\inf_{q\in\mathbb{R}^n}\big\{\lambda_{\min}(J(q)^\top\mathcal{Q}(q)^\top \mathcal{Q}(q) J(q))\big\}>0$.
\end{definition}

to be satisfied by the stacked Jacobian $J$ and the orientation map $\mathcal{Q}$, i.e.\ by the kinematic layout of the robot. This can be seen by means of the following counterexample: The inertial parameters $m_l=1$, $h_l=(\begin{smallmatrix}
        0&0&0
\end{smallmatrix})^{ \top}$ and $I_l=\mathrm{I}_3$ for $l\in\{1,\ldots,N\}$ are constant and physically consistent (the latter is readily verified via Theorem~\ref{theorem:PhysicalConsistency} since $f(\Phi_l)=\mathrm{diag} (
        \tfrac{1}{2}\mathrm{I}_3, 1)\succ 0$) and they lead to $\mathcal{Z}(\Theta)=\mathrm{I}_{6N}$, cf.~\eqref{DefMatrixZ},~\eqref{DefBlockmatrixZ}. According to~\eqref{MassMatrix}, the mass matrix with these inertial parameters is $M(q,\Theta)=J(q)^\top\mathcal{Q}(q)^\top \mathcal{Q}(q) J(q)$ such that loss of normality of $\mathcal{Q}(\cdot)J(\cdot)$ would cause $\inf_{q\in\mathbb{R}^n}(\lambda_{\min}(M(q,\Theta)))=0$, i.e.\ causing the greatest lower uniform bound of the mass matrix to lose its desirable positivity despite the fact that the inertial parameters are chosen physically consistent and constant.

Moreover, it can be shown for constant, physically consistent inertial parameters that a positive lower uniform bound of the mass matrix exists if and only if $\mathcal{Q}(\cdot)J(\cdot)$ is normal (see Lemma~\ref{lemma:NormalJacobian} from below) and that a finite upper uniform bound of the mass matrix exists if and only if the robotic manipulator meets the specifications listed in \cite{ghorbel_uniform_1998}.

However, when dealing with time-dependent inertial parameters, then their physical consistency at all times, compliance with the specifications listed in \cite{ghorbel_uniform_1998} and normality of $\mathcal{Q}(\cdot)J(\cdot)$ are not anymore sufficient to guarantee the existence of finite, positive uniform bounds of the mass matrix, as is made evident by yet another counterexample: Consider a robotic manipulator with $N=1$ describing the translational movement of a solid uniform sphere with mass $m_1(t)\in\mathbb{R}$ and radius $r>0$, i.e.\ $\rho_1(x_1,t)=\tfrac{3}{4\pi r^3} m_1(t)$ for $\|x_1\|\leq r$ and $\rho_1(x_1,t)=0$ elsewhere, along the $x$-axis of the $0$th Newtonian frame. This results in the stacked Jacobian $J(q)=(\begin{smallmatrix}
        1&0&0&0&0&0
    \end{smallmatrix})^{ \top}$, ensures normality of $\mathcal{Q}(\cdot)J(\cdot)$ since $\mathcal{Q}(q)J(q)=((\begin{smallmatrix}
        1&0&0
    \end{smallmatrix})\mathcal{R}(F_\mathrm{\phi_1}(q))^\top,(\begin{smallmatrix}
        0&0&0
    \end{smallmatrix}))^\top$ and $\mathcal{R}(F_\mathrm{\phi_1}(q))$ is orthogonal and it renders the mass matrix $M(q,\Theta)= m_1$, cf.~\eqref{DefMatrixZ},~\eqref{MassMatrix},~\eqref{DefBlockmatrixZ}. The other inertial parameters of the sphere besides $m_1(t)$ are $h_1(t)=(\begin{smallmatrix}
        0&0&0
    \end{smallmatrix})^{ \top}$ and $I_1(t)=\tfrac{2}{5}m_1(t)r^2\mathrm{I}_3$. This leads to $f(\Phi_1(t))=m_1(t)\,\mathrm{diag} (
        \tfrac{1}{5}r^2\mathrm{I}_3, 1
    )$ and thus shows that these inertial parameters at all times are physically consistent if and only if $m_1(t)>0$ for all $t\geq 0$. Thereby, physical consistency at all times neither prevents the greatest lower uniform bound of $M$ from being nonpositive when $\lim_{t\to \infty}m_1(t)=0$ nor the lowest upper uniform bound of $M$ from being infinity when $\lim_{t\to \infty}m_1(t)=\infty$ even though this robotic manipulator meets the requirements listed in \cite{ghorbel_uniform_1998}.

To circumvent difficulties that arise when the inertial parameters approach physical inconsistency or diverge, we introduce some additional properties that rule out such pathological behavior:

\begin{definition}\label{definition:PropertiesInertialParameters}
     The inertial parameters are \textit{uniformly physically consistent} if $\inf_{t\geq 0}\big\{\lambda_{\min}(f(\Phi_l(t)))\big\}>0$ for all $l\in\{1,\ldots,N\}$ and they are \textit{upper bounded} if $\sup_{t\geq 0}\big\{\lambda_{\max}(f(\Phi_l(t)))\big\}<\infty$ for all $l\in\{1,\ldots,N\}$.
\end{definition}   

Equipped with the Definitions~\ref{definition:NormalJacobian},~\ref{definition:PropertiesInertialParameters}, the standing Assumption~\ref{Assumption:IndependentJacobian} allows us---provided some technical presuppositions hold---to derive finite, positive uniform bounds for the mass matrix from upper boundedness and {uniform} physical consistency of time-dependent inertial parameters. Specifically, we investigate the existence of a positive lower uniform bound of the mass matrix and then connect the presuppositions in \cite{ghorbel_uniform_1998} for the existence of a finite upper uniform bound of the mass matrix for constant inertial parameters to the time-dependent case. Afterward, these results are summarized and discussed in a unified~fashion.

\paragraph{Positive lower uniform bound}

\begin{theorem}\label{theorem:LowerBoundedMassMatrix}
    If $\mathcal{Q}(\cdot)J(\cdot)$ is normal and the inertial parameters are uniformly physically consistent, then the mass matrix~\eqref{MassMatrix} is uniformly bounded by a positive lower bound, specifically then $M(q,\Theta(t))\succeq  \alpha_1 \mathrm{I}_n$ for all $t\geq 0$, $q\in\mathbb{R}^n$ with any positive constant
    \begin{align}
    \label{LowerBoundedness}
        \alpha_1< \underbrace{\min_{l\in\{1,\ldots,N\}}\Big\{\inf_{t\geq 0}\big\{\lambda_{\min}(f(\Phi_l(t)))\big\}\Big\}}_{\text{positive due to uniform physical consistency}}\, \underbrace{\inf_{q\in\mathbb{R}^n}\big\{\lambda_{\min}(J(q)^\top \mathcal{Q}(q)^\top\mathcal{Q}(q) J(q))\big\}}_{\text{positive due to normality}}.
    \end{align}
\end{theorem}
\begin{proof}
    The statements in this proof are valid for all $l\in\{1,\ldots,N\}$, $t\geq 0$. Suppose, the inertial parameters are uniformly physically consistent at all times, i.e.\ we find constants $\xi_l>0$ with $f(\Phi_l(t))\succ \xi_l \mathrm{I}_4$. Thereby, $\Big(\begin{smallmatrix}\Sigma_l(t)-\xi_l\mathrm{I}_3&h_l(t)\\h_l(t)^\top&m_l(t)-\xi_l\end{smallmatrix}\Big)\succ 0$ such that we receive $m_l(t)-\xi_l>0$ and $\Sigma_l(t)-\xi_l\mathrm{I}_3-(m_l(t)-\xi_l)^{-1}h_l(t)h_l(t)^\top\succ 0$ by application of Schur's complement and thus arrive at $\mathrm{tr}(\Sigma_l(t)-\xi_l\mathrm{I}_3-(m_l(t)-\xi_l)^{-1}h_l(t)h_l(t)^\top)\mathrm{I}_3\succ \Sigma_l(t)-\xi_l\mathrm{I}_3-(m_l(t)-\xi_l)^{-1}h_l(t)h_l(t)^\top$. Reorganization gives $\mathrm{tr}(\Sigma_l(t))\mathrm{I}_3-\Sigma_l(t)-(m_l(t)-\xi_l)^{-1}(\mathrm{tr}(h_l(t)h_l(t)^\top)\mathrm{I}_3-h_l(t)h_l(t)^\top)\succ \mathrm{tr}(\xi_l\mathrm{I}_3)\mathrm{I}_3-\xi_l\mathrm{I}_3\succ \xi_l\mathrm{I}_3$. Plugging in the inverse transformation from Theorem~\ref{theorem:PhysicalConsistency} while using the equality $\mathcal{S}(h_l(t))^\top \mathcal{S}(h_l(t))=\mathrm{tr}(h_l(t)h_l(t)^\top)\mathrm{I}_3-h_l(t)h_l(t)^\top$, which itself is readily verified, leads to $I_l(t)-(m_l(t)-\xi_l)^{-1}\mathcal{S}(h_l(t))^\top \mathcal{S}(h_l(t))\succ \xi_l \mathrm{I}_3$, i.e.\ $I_l(t)-\xi_l \mathrm{I}_3-(-\mathcal{S}(h_l(t))^\top) ((m_l(t)-\xi_l)\mathrm{I}_3)^{-1}(-\mathcal{S}(h_l(t)))\succ 0$. Application of Schur's complement with respect to the last inequality and $(m_l(t)-\xi_l)\mathrm{I}_3\succ 0$ reveals $\Big(\begin{smallmatrix}(m_l(t)-\xi_l)\mathrm{I}_3&-\mathcal{S}(h_l(t))\\-\mathcal{S}(h_l(t))^\top&I_l(t)-\xi_l\mathrm{I}_3\end{smallmatrix}\Big)\succ 0$, i.e.\ $\lambda_{\min}(Z(\Phi_l(t)))\geq \xi_l$, cf.\ the definition of the matrix $Z(\Phi_l(t))$ in~\eqref{DefMatrixZ}. By means of~\eqref{MassMatrix},~\eqref{DefBlockmatrixZ}, we arrive at $\dot{q}^\top M(q,\Theta(t))\dot{q}=\dot{q}^\top J(q)^\top\mathcal{Q}(q)^\top \mathrm{diag}(Z(\Phi_1(t)), \ldots, Z(\Phi_N(t))) \mathcal{Q}(q)J(q)\dot{q}\geq \beta \dot{q}^\top J(q)^\top\mathcal{Q}(q)^\top \mathcal{Q}(q) J(q)\dot{q}\geq \beta\gamma\|\dot{q}\|^2$ for all $q,\dot{q}\in\mathbb{R}^n$ with the constants $\beta=\min_{l\in\{1,\ldots,N\}}\big\{\xi_l\big\}>0$, $\gamma=\inf_{q\in\mathbb{R}^n}\big\{\lambda_{\min}(J(q)^\top\mathcal{Q}(q)^\top \mathcal{Q}(q) J(q))\big\}$. Accordingly, $\lambda_{\min}(M(q,\Theta(t)))\geq \beta\gamma=\alpha_1$, where the constant $\alpha_1$ is positive if $\mathcal{Q}(\cdot)J(\cdot)$ is normal, i.e.\ $\gamma>0$, in addition to the presupposed positivity of $\xi_l$ and $\beta$. Moreover, we may rewrite $\alpha_1=\gamma\min_{l\in\{1,\ldots,N\}}\{\xi_l\}$ such that~\eqref{LowerBoundedness} follows by recalling that $\xi_l$ can be chosen arbitrarily in the open interval $\big(0,\inf_{t\geq 0}\{\lambda_{\min}(f(\Phi_l(t)))\}\big)$.
\end{proof}

The results in Theorem~\ref{theorem:LowerBoundedMassMatrix} lead to the next statement that provides a characterization of normality of $\mathcal{Q}(\cdot)J(\cdot)$ by means of constant, physically consistent inertial parameters:
\begin{lemma}\label{lemma:NormalJacobian}
	$\mathcal{Q}(\cdot)J(\cdot)$ is normal if and only if the greatest lower uniform bound of the mass matrix~\eqref{MassMatrix} is positive, i.e.\ $\inf_{q\in\mathbb{R}^n}\big\{\lambda_{\min}(M(q,\Theta))\big\}>0$, whenever the inertial parameters with $f(\Phi_l)\succ  0$ for all $l\in\{1,\ldots,N\}$ are physically consistent and constant.
\end{lemma}
\begin{proof}
    The sufficient direction represents a special case of Theorem~\ref{theorem:LowerBoundedMassMatrix}. The necessary direction follows when considering the constant, physically consistent inertial parameters from the first counterexample with $\mathcal{Z}(\Theta)=\mathrm{I}_{6N}$ because they cause $\inf_{q\in\mathbb{R}^n}\{\lambda_{\min}(J(q)^\top\mathcal{Q}(q)^\top \mathcal{Q}(q) J(q))\}=\inf_{q\in\mathbb{R}^n}\{\lambda_{\min}(M(q,\Theta))\}$, cf.~\eqref{MassMatrix}.
\end{proof}

\paragraph{Finite upper uniform bound}\begin{assumption}\label{Assumption:UniformUpperBound}
    If the inertial parameters are physically consistent and constant, i.e.\ $f(\Phi_l)\succ 0$ for all $l\in\{1,\ldots,N\}$, then $\sup_{q\in\mathbb{R}^n}\big\{\lambda_{\max}(M(q,\Theta))\big\}<\infty$. 
\end{assumption}
In \cite{ghorbel_uniform_1998}, Fathi Gorbel et al.\ exactly characterize open kinematic chains that possess the desirable property $\sup_{q\in\mathbb{R}^n}\big\{\lambda_{\max}(M(q,\Theta))\big\}<\infty$, provided that the analysis is restricted to constant, physically consistent inertial parameters. More precisely, the joints must appear in a certain order, depending on whether they are revolute or prismatic and the rotations of the rigid bodies must maintain a certain orientation relative to each other for all $q\in\mathbb{R}^n$. This means that the specifications listed in \cite{ghorbel_uniform_1998} put constraints on the poses of the joints relative to each other. Therefore---in light of the standing Assumption~\ref{Assumption:IndependentJacobian} and provided that the origins of the body-fixed frames are located at the joints---we can evaluate the validity of these constraints by means of information contained in the forward kinematic map $F$ and the underlying orientation map $\mathcal{Q}$ independently of the mass densities and thus independently of the inertial parameters. Hence, if all poses of the joints that are possible during operation of the robotic manipulator for one arbitrary set of constant, physically consistent inertial parameters meet the specifications listed in \cite{ghorbel_uniform_1998}, then the Assumption~\ref{Assumption:UniformUpperBound} is fulfilled. This simple way of checking the validity of Assumption~\ref{Assumption:UniformUpperBound} opens up the possibility to assess the existence of a finite upper uniform bound of the mass matrix even in the more involved case of time-dependent inertial parameters by means of the next result:

\begin{theorem}\label{theorem:UpperBoundedMassMatrix}
    Suppose, Assumption~\ref{Assumption:UniformUpperBound} is fulfilled. If the inertial parameters are physically consistent at all times and upper bounded, then the mass matrix~\eqref{MassMatrix} is uniformly bounded through a finite upper bound, specifically then there exists a constant $0\leq\alpha_2<\infty$ with $0\preceq  M(q,\Theta(t))\preceq \alpha_2 \mathrm{I}_n$ for all $q\in\mathbb{R}^n$, $t\geq 0$ that satisfies
    \begin{align}
    \label{UpperBoundedness}
        \alpha_2\leq 2\underbrace{\max_{l\in\{1,\ldots,N\}}\Big\{\sup_{t\geq 0}\big\{\lambda_{\max}(f(\Phi_l(t)))\big\}\Big\}}_{\text{finite due to upper boundedness}}\,\underbrace{\sup_{q\in\mathbb{R}^n}\big\{\sigma_{\max}(\mathcal{Q}(q)J(q))^2\big\}}_{\text{finite as per Assumption~\ref{Assumption:UniformUpperBound}, cf.\ Lemma~\ref{lemma:UpperBoundedJacobian}}}.
    \end{align}
\end{theorem}
\begin{proof}
    Consider inertial parameters with $f(\Phi_l(t))\succ 0$ and $\sup_{t\geq 0}\big\{\lambda_{\max}(f(\Phi_l(t)))\big\}<\infty$ for all $l\in\{1,\ldots,N\}$, $t\geq 0$, i.e.\ they are physically consistent at all times and upper bounded. Since Assumption~\ref{Assumption:UniformUpperBound} is fulfilled and the inertial parameters at $t=0$ are physically consistent, i.e.\ $f(\Phi_l(0))\succ 0$ for all $l\in\{1,\ldots,N\}$, we find a constant $0\leq\zeta< \infty$ with $\zeta\|\dot{q}\|^2\geq \dot{q}^\top M(q,\Theta(0))\dot{q}=\dot{q}^\top J(q)^\top\mathcal{Q}(q)^\top \mathcal{Z}(\Theta(0))\cdot\allowbreak\mathcal{Q}(q) J(q)\dot{q}\geq \lambda_{\min}(\mathcal{Z}(\Theta(0)))\dot{q}^\top J(q)^\top\mathcal{Q}(q)^\top\mathcal{Q}(q) J(q)\dot{q}$ for all $q,\dot{q}\in\mathbb{R}^n$, cf.~\eqref{MassMatrix}. Following the argumentation in the proof of Theorem~\ref{theorem:LowerBoundedMassMatrix}, physical consistency of the inertial parameters at $t=0$ causes $\lambda_{\min}(Z(\Phi_l(0)))>0$ for all $l\in\{1,\ldots,N\}$ and consequently it leads to $\lambda_{\min}(\mathcal{Z}(\Theta(0)))>0$ for the corresponding block structure. We arrive at $\dot{q}^\top J(q)^\top \mathcal{Q}(q)^\top \mathcal{Q}(q) J(q)\dot{q}\leq \frac{\zeta}{\lambda_{\min}(\mathcal{Z}(\Theta(0)))}\|\dot{q}\|^2$ for all $q,\dot{q}\in\mathbb{R}^n$ and thus obtain $\sup_{q\in\mathbb{R}^n}\big\{\lambda_{\max}(J(q)^\top \mathcal{Q}(q)^\top \mathcal{Q}(q) J(q))\big\}\leq\delta$ with the constant $\delta=\frac{\zeta}{\lambda_{\min}(\mathcal{Z}(\Theta(0)))}\in[0,\infty)$. The subsequent argumentation aims to infer upper boundedness of $\mathcal{Z}(\Theta(\cdot))$ from the presupposed upper boundedness of the inertial parameters. The remaining statements in this proof are valid for all $l\in\{1,\ldots,N\}$, $t\geq 0$. Due to $\lambda_{\max}(f(\Phi_l(t)))<\infty$, we find constants $\zeta_l\in[0,\infty)$ with $f(\Phi_l(t))\prec \zeta_l\mathrm{I}_4$, i.e.\ $\Big(\begin{smallmatrix}
        \zeta_l\mathrm{I}_3-\Sigma_l(t)&-h_l(t)\\-h_l(t)^\top&\zeta_l-m_l(t)
    \end{smallmatrix}\Big)\succ 0$. According to Schur's complement, this implies $\zeta_l-m_l(t)>0$ and $\zeta_l\mathrm{I}_3-\Sigma_l(t)-(\zeta_l-m_l(t))^{-1}h_l(t)h_l(t)^\top\succ 0$, where the latter inequality leads to $\mathrm{tr}(\zeta_l\mathrm{I}_3-\Sigma_l(t)-(\zeta_l-m_l(t))^{-1}h_l(t)h_l(t)^\top)\mathrm{I}_3\succ \zeta_l\mathrm{I}_3-\Sigma_l(t)-(\zeta_l-m_l(t))^{-1}h_l(t)h_l(t)^\top$. Reordering shows $\Sigma_l(t)-\mathrm{tr}(\Sigma_l(t))\mathrm{I}_3-(\zeta_l-m_l(t))^{-1}(\mathrm{tr}(h_l(t)h_l(t)^\top)\mathrm{I}_3-h_l(t)h_l(t)^\top)\succ \zeta_l\mathrm{I}_3-\mathrm{tr}(\zeta_l\mathrm{I}_3)\mathrm{I}_3=-2\zeta_l\mathrm{I}_3$, which, by taking into account the inverse transformation in Theorem~\ref{theorem:PhysicalConsistency} and $\mathcal{S}(h_l(t))^\top \mathcal{S}(h_l(t))=\mathrm{tr}(h_l(t)h_l(t)^\top)\mathrm{I}_3-h_l(t)h_l(t)^\top$, rewrites as $I_l(t)+(\zeta_l-m_l(t))^{-1}\mathcal{S}(h_l(t))^\top \mathcal{S}(h_l(t))\prec 2\zeta_l\mathrm{I}_3$. Due to $\zeta_l\geq 0$, the more conservative inequality $I_l(t)+(2\zeta_l-m_l(t))^{-1}\mathcal{S}(h_l(t))^\top \mathcal{S}(h_l(t))\prec 2\zeta_l\mathrm{I}_3$ is also valid, i.e.\ we obtain $2\zeta_l\mathrm{I}_3-I_l(t)-\mathcal{S}(h_l(t))^\top((2\zeta_l-m_l(t))\mathrm{I}_3)^{-1}\mathcal{S}(h_l(t))\succ 0$. Applying Schur's complement with respect to the last inequality and the consequence $(2\zeta_l-m_l(t))\mathrm{I}_3\succ 0$ of $\zeta_l\geq 0$ and $\zeta_l-m_l(t)>0$ yields $\big(\begin{smallmatrix}
        (2\zeta_l-m_l(t))\mathrm{I}_3&\mathcal{S}(h_l(t))\\\mathcal{S}(h_l(t))^\top&2\zeta_l\mathrm{I}_3-I_l
    \end{smallmatrix}\big)\succ 0$. This leads to $Z(\Phi_l(t))\prec 2\zeta_l\mathrm{I}_6$, cf.\ the definition in~\eqref{DefMatrixZ}, i.e.\ we have verified upper boundedness of $\mathcal{Z}(\Theta(\cdot))$ with $\varepsilon=\sup_{t\geq 0}\big\{\lambda_{\max}(\mathcal{Z}(\Theta(t)))\big\}=\sup_{t\geq 0,\,l\in\{1,\ldots,N\}}\big\{\lambda_{\max}(Z(\Phi_l(t)))\big\}\leq \max_{l\in\{1,\ldots,N\}}\big\{2\zeta_l\big\}\in[0,\infty)$. Hence, the structure of~\eqref{MassMatrix} shows $\dot{q}^\top\! M(q,\Theta(t))\dot{q}=\dot{q}^\top\! J(q)^\top\!\mathcal{Q}(q)^\top \mathcal{Z}(\Theta(t)) \mathcal{Q}(q) J(q)\dot{q}\leq\varepsilon \dot{q}^\top J(q)^\top \mathcal{Q}(q)^\top \mathcal{Q}(q) J(q)\dot{q} \leq \varepsilon \sigma_{\max}\big\{\mathcal{Q}(q)J(q)\}^2\|\dot{q}\|^2\leq \alpha_2 \|\dot{q}\|^2$ for all $q,\dot{q}\in\mathbb{R}^n$ with $\alpha_2=\varepsilon\sup_{q\in\mathbb{R}^n}\big\{\sigma_{\max}\{\mathcal{Q}(q)J(q)\}^2\big\}\leq \varepsilon\delta \in[0,\infty)$, revealing the constant $\alpha_2$ to be nonnegative and finite, as desired. Furthermore, $\zeta_l$ can be chosen in compliance with its initial definition as $\zeta_l=\varsigma \sup_{t\geq 0}\{\lambda_{\max}(f(\Phi_l(t)))$ with any $\varsigma>1$. Thus, by recalling the structure of $\alpha_2$ and $\varepsilon$, we obtain $\alpha_2\leq 2 \max_{l\in\{1,\ldots,N\}}\{\zeta_l\}\sup_{q\in\mathbb{R}^n}\big\{\sigma_{\max}\{\mathcal{Q}(q)J(q)\}^2\big\}=2 \varsigma \max_{l\in\{1,\ldots,N\}}\{\sup_{t\geq 0}\{\lambda_{\max}(f(\Phi_l(t)))\}\sup_{q\in\mathbb{R}^n}\big\{\sigma_{\max}\{\mathcal{Q}(q)J(q)\}^2\big\}$ for any $\varsigma>1$, which implies~\eqref{UpperBoundedness} when choosing $\varsigma$ arbitrarily close to one.
\end{proof}

The next lemma shows how to check the validity of Assumption~\ref{Assumption:UniformUpperBound}, i.e.\ the condition to obtain a finite upper uniform bound for the mass matrix by means of Theorem~\ref{theorem:UpperBoundedMassMatrix}, without taking a detour via \cite{ghorbel_uniform_1998} but instead directly using information carried in the forward kinematic map $F$, specifically the stacked Jacobian $J$, together with the orientation map $\mathcal{Q}$:

\begin{lemma}\label{lemma:UpperBoundedJacobian}
    $\mathcal{Q}(\cdot)J(\cdot)$ is bounded with $\sup_{q\in\mathbb{R}^n}\big\{\sigma_{\max}(\mathcal{Q}(q)J(q))\big\}<\infty$ if and only if the mass matrix~\eqref{MassMatrix} fulfills Assumption~\ref{Assumption:UniformUpperBound}. 
\end{lemma}
\begin{proof}
    First, we verify the sufficient direction of the desired equivalence. The structure of the mass matrix in~\eqref{MassMatrix} shows $\dot{q}^\top M(q,\Theta)\dot{q}\leq \sigma_{\max}(\mathcal{Q}(q)J(q))^2\lambda_{\max}(\mathcal{Z}(\Theta))\|\dot{q}\|^2$ for all $q,\dot{q}\in\mathbb{R}^n$ and arbitrary inertial parameters $\Theta\in\mathbb{R}^{10N}$. This results in $\lambda_{\max}(M(q,\Theta))\leq\sigma_{\max}(\mathcal{Q}(q)J(q))^2\lambda_{\max}(\mathcal{Z}(\Theta))$.  Since physically consistent, constant inertial parameters trivially are upper bounded, they admit $\lambda_{\max}(\mathcal{Z}(\Theta))\in(0,\infty)$, cf.\ the proof of Theorem~\ref{theorem:UpperBoundedMassMatrix}, and thus, we receive $\sup_{q\in\mathbb{R}^n}\big\{\lambda_{\max}(M(q,\Theta))\}<\infty$ for such inertial parameters if $\sup_{q\in\mathbb{R}^n}\big\{\sigma_{\max}(\mathcal{Q}(q)J(q))\big\}<\infty$. Next, in order to show necessity in the desired equivalence, consider rigid bodies with the inertial parameters $m_l=1$, $h_l=(\begin{smallmatrix}
        0&0&0
    \end{smallmatrix})^{ \top}$ and $I_l=\mathrm{I}_3$ for $l\in\{1,\ldots,N\}$ such that $\mathcal{Z}(\Theta)=\mathrm{I}_{6N}$. Then, the structure of~\eqref{MassMatrix} leads to $\sup_{q\in\mathbb{R}^n}\big\{\lambda_{\max}(M(q,\Theta))\big\}=\sup_{q\in\mathbb{R}^n}\big\{\lambda_{\max}(J(q)^\top\mathcal{Q}(q)^\top \mathcal{Q}(q) J(q))\big\}=\sup_{q\in\mathbb{R}^n}\big\{\sigma_{\max}(\mathcal{Q}(q)J(q))\big\}^2$. Accordingly, the physical consistency of these specific inertial parameters, which is confirmed by means of Theorem~\ref{theorem:PhysicalConsistency} in view of $f(\Phi_l)=\mathrm{diag} (
        \tfrac{1}{2}\mathrm{I}_3, 1
    )\succ 0$, shows that validity of Assumption~\ref{Assumption:UniformUpperBound} implies $\sup_{q\in\mathbb{R}^n}\big\{\sigma_{\max}(\mathcal{Q}(q)J(q))\big\}<\infty$.
\end{proof}

\paragraph{Unifying discussion of mass-matrix boundedness} Irrespective of the way one evaluates the validity
of Assumption~\ref{Assumption:UniformUpperBound}, the summary of the findings on the existence of finite, positive uniform bounds of the mass matrix, i.e.\ the conjunction of Theorem~\ref{theorem:LowerBoundedMassMatrix} and Theorem~\ref{theorem:UpperBoundedMassMatrix}, presents itself as:
\begin{corollary}\label{corollary:UniformBounds}
    Suppose, $\mathcal{Q}(\cdot)J(\cdot)$ is normal and Assumption~\ref{Assumption:UniformUpperBound} is fulfilled. If the inertial parameters are upper bounded and uniformly physically consistent, then the mass matrix~\eqref{MassMatrix} is uniformly bounded with $\alpha_1 \mathrm{I}_n\preceq M(q,\Theta(t))\preceq\alpha_2\mathrm{I}_n$ for all $q\in\mathbb{R}^n$, $t\geq 0$, where $0<\alpha_1\leq\alpha_2<\infty$ are constants.
\end{corollary}
Further, by conjunction of Lemma~\ref{lemma:NormalJacobian} and Lemma~\ref{lemma:UpperBoundedJacobian}, we obtain the subsequent result that provides contrast to the time-dependent case described in Corollary~\ref{corollary:UniformBounds} when the inertial parameters are constant:
\begin{corollary}\label{corollary:UniformBoundsForConstantParameters}
	$\mathcal{Q}(\cdot)J(\cdot)$ is normal and bounded, i.e.\ it admits $\beta_1\mathrm{I}_n\preceq J(q)^\top\mathcal{Q}(q)^\top \mathcal{Q}(q) J(q)\preceq \beta_2\mathrm{I}_n$ for all $q\in\mathbb{R}^n$ with some constants $0<\beta_1\leq\beta_2<\infty$, if and only if there exist finite, positive uniform bounds $0<\alpha_1\leq\alpha_2<\infty$ of the mass matrix~\eqref{MassMatrix} with $\alpha_1 \mathrm{I}_n\preceq M(q,\Theta)\preceq\alpha_2\mathrm{I}_n$ for all $q\in\mathbb{R}^n$ whenever the inertial parameters with $f(\Phi_l)\succ 0$ for all $l\in\{1,\ldots,N\}$ are constant and physically consistent.
\end{corollary}
Notably, Corollary~\ref{corollary:UniformBoundsForConstantParameters} is an equivalence. This means that the required normality of $\mathcal{Q}(\cdot)J(\cdot)$ and fulfillment of Assumption~\ref{Assumption:UniformUpperBound}---i.e.\ additional boundedness of $\mathcal{Q}(\cdot)J(\cdot)$, see Lemma~\ref{lemma:UpperBoundedJacobian}---for existence of finite, positive uniform bounds of the mass matrix in Corollary~\ref{corollary:UniformBounds} are not only sufficient but also necessary for validity of an assertion commonly made implicitly by contributions on robot control, e.g.\ in~\cite{slotine_adaptive_1987,wu_adaptive_2022}: that finite, positive uniform bounds of the mass matrix will exist simply because the inertial parameters correspond to physically meaningful rigid bodies, i.e.\ because they are physically consistent and constant. Thus, if one is willing to assert existence of finite, positive uniform bounds of the mass matrix based purely on rigidity of and physical meaning in the bodies that build a robotic manipulator, then Corollary~\ref{corollary:UniformBounds} may be invoked to discuss the more involved case of time-dependent inertial parameters without need to impose any further restriction on $\mathcal{Q}(\cdot)J(\cdot)$ representing the kinematic layout of such a robot. Moreover, on the contrary, should a given robotic manipulator break this common assertion because its kinematic layout lacks either normality or boundedness of the corresponding $\mathcal{Q}(\cdot)J(\cdot)$, then either Theorem~\ref{theorem:UpperBoundedMassMatrix} or Theorem~\ref{theorem:LowerBoundedMassMatrix} still provides the means to assess the opposing uniform bound of the mass matrix by judging either upper boundedness or uniformity of physically consistent inertial parameters, respectively. Thereby, the only case for which none of our results on the boundedness of the mass matrix are applicable is the degenerate one where the common assertion from above breaks because the kinematic layout of a given robotic manipulator is described by a non-normal $\mathcal{Q}(\cdot)J(\cdot)$ that is not bounded.

\subsection{Boundedness of other components in the generalized 
robotics equation}\label{section:OtherComponents}
The estimation regimes of the adaptation schemes in  \cite{craig_adaptive_1987,slotine_adaptive_1987,johansson_adaptive_1990,lee_natural_2018} are derived under the presumption that the parameters $\Theta$ might be unknown but remain constant, i.e.\ assuming absence of the effect $M(q,\dot{\Theta})\dot{q}$ from the dynamical behavior as  represented by the generalized robotics equation~\eqref{RoboticEquation}. Therefore, testing the robustness of these algorithms should involve a violation of this assumption and thus requires bounding the term $M(q,\dot{\Theta})\dot{q}$ in cases where only existence of a finite bound of the rate of change of the unknown parameters can be guaranteed. Now, the fact that here in this work the values in $\Theta$ are inertial parameters allows exploiting the resulting structure~\eqref{MassMatrix} of $M(q,\dot{\Theta})$ and thereby enables us to deliver such a boundedness result down below in Theorem~\ref{theorem:BoundedRateOfChange}. Therein, we choose to associate the bound of the rate of change of the inertial parameters $\Phi_l$ with $\sup_{t\geq 0}\big\{\sigma_{\max}(f(\dot{\Phi}_l(t)))\big\}$ because the maximal singular value offers indifference toward the direction of the changes occurring in $\Phi_l$ over time and we deem these directions to be irrelevant for the intended robustness testing.
\begin{theorem}\label{theorem:BoundedRateOfChange}
    Let Assumption~\ref{Assumption:UniformUpperBound} apply and the mass matrix be structured as in~\eqref{MassMatrix}. If the time derivatives of the inertial parameters are bounded with $\sup_{t\geq 0}\big\{\sigma_{\max}(f(\dot{\Phi}_l(t)))\big\}<\infty$ for all $l\in\{1,\ldots,N\}$, then the matrix $M(q,\dot{\Theta})$ is uniformly bounded through a finite bound, specifically the finite, nonnegative constant
    \begin{align}
    \label{BoundedEffectofVaryingParam}
    \chi = \sqrt{69N} \underbrace{\max_{l\in\{1,\ldots,N\}}\Big\{\sup_{t\geq 0}\big\{\sigma_{\max}(f(\dot{\Phi}_l(t)))\big\}\Big\}}_{\text{finite due to bounded rates of change}} \,\underbrace{\sup_{q\in\mathbb{R}^n}\big\{\sigma_{\max}(\mathcal{Q}(q)J(q))^2\big\}}_{\text{finite as per Assumption~\ref{Assumption:UniformUpperBound}, cf.\ Lemma~\ref{lemma:UpperBoundedJacobian}}}
    \end{align}
    then admits $\sigma_{\max}(M(q,\dot{\Theta}(t)))\leq \chi$ for all $q\in\mathbb{R}^n$, $t\geq 0$.
\end{theorem}
\begin{proof}
    The statements in this proof hold for all $l\in\{1,\ldots,N\}$, $t\geq 0$. Consider inertial parameters whose time derivatives admit $\mu=\max_{l\in\{1,\ldots,N\}}\big\{\sup_{t\geq 0}\big\{\sigma_{\max}(f(\dot{\Phi}_l(t)))\big\}\big\}<\infty$. Accordingly, we have $\sigma_{\max}(f(\dot{\Phi}_l(t)))\leq\mu$ and thereby, the absolute values of the entries in $f(\dot{\Phi}_l(t))$ are bounded by $\mu$,~i.e.\ $ |\dot{m}_l(t) |, |\dot{h}_{\mathrm{x},l}(t) |, |\dot{h}_{\mathrm{y},l}(t) |, |\dot{h}_{\mathrm{z},l}(t) |\leq \mu$ and $ |\mathrm{e}_{i,3}^\top \dot{\Sigma}_l(t)\mathrm{e}_{j,3} |= |\mathrm{tr}(\dot{\Sigma}_l(t))\mathrm{e}_{i,3}^\top \mathrm{e}_{j,3}-\mathrm{e}_{i,3}^\top \dot{I}_l(t)\mathrm{e}_{j,3} |\leq\mu$ for all $i,j\in\{1,2,3\}$ when taking into account the consequence $\mathrm{tr}(\dot{I}_l(t))=2\mathrm{tr}(\dot{\Sigma}_l(t))$ of the inverse transformation in Theorem~\ref{theorem:PhysicalConsistency}. Therein, $\mathrm{e}_{i,\mathscr{n}}\in\mathbb{R}^{\mathscr{n}}$ is the $i$th standard unit vector. For $i<j$, we receive $ |\dot{I}_{ij,l}(t) |\leq\mu$ for the entries in $\dot{I}_l(t)$ offside its main diagonal. For $i=j$, we obtain $ |\mathrm{tr}(\dot{\Sigma}_l(t))-\dot{I}_{ii,l}(t) |\leq \mu$ such that $ |\mathrm{tr}(\dot{\Sigma}_l(t)) |\leq 3\mu$ necessitates the diagonal entries of  $\dot{I}_l(t)$ to satisfy $ |\dot{I}_{ii,l}(t) |\leq 4\mu$. Thereby, we obtain $\sigma_{\max}(\mathcal{Z}(\dot{\Theta}(t)))=(\lambda_{\max}(\mathcal{Z}(\dot{\Theta}(t))^\top \mathcal{Z}(\dot{\Theta}(t))))^{1/2}\leq (\mathrm{tr}(\mathcal{Z}(\dot{\Theta}(t))^\top \mathcal{Z}(\dot{\Theta}(t))))^{1/2}=(\sum_{{j}=1}^{6N}(\mathcal{Z}(\dot{\Theta}(t))\mathrm{e}_{{j},6N})^\top \mathcal{Z}(\dot{\Theta}(t))\mathrm{e}_{{j},6N})^{1/2}=(\sum_{{j}=1}^{6N}\sum_{{i}=1}^{6N}(\mathrm{e}_{{i},6N}^\top\mathcal{Z}(\dot{\Theta}(t))\mathrm{e}_{{j},6N})^2)^{1/2}= (\sum_{l=1}^N\sum_{{j}=1}^{6}\sum_{{i}=1}^{6}|\mathrm{e}_{{i},6}^\top{Z}(\dot{\Phi}_l(t))\mathrm{e}_{{j},6}|^2)^{1/2}\leq (N\cdot(2\cdot 6+3+6+ 3\cdot 4^2)\mu^2)^{1/2}=\sqrt{69N}\mu$ in view of the block-diagonal structure~\eqref{DefBlockmatrixZ} of $\mathcal{Z}(\dot{\Theta}(t))$ with the matrices $Z(\dot{\Phi}_l(t))$ on its main diagonal, that as per~\eqref{DefMatrixZ} have $2\cdot 6$ entries from the elements of $\dot{h}_l(t)$, $3$ duplicate entries of $\dot{m}_l(t)$, $6$ entries from the off-diagonal elements of $\dot{I}_{l}(t)$ and the $3$ diagonal entries of $\dot{I}_{l}(t)$, respectively; all of which sums up to the upper bound $\sqrt{69 N}\mu$, as stated. Further, we get $ \dot{q}^\top M(q,\dot{\Theta}(t))^\top M(q,\dot{\Theta}(t)) \dot{q}=\dot{q}^\top J(q)^\top\mathcal{Q}(q)^\top  \mathcal{Z}(\dot{\Theta}(t))^\top \mathcal{Q}(q) J(q) J(q)^\top \mathcal{Q}(q)^\top \mathcal{Z}(\dot{\Theta}(t)) \mathcal{Q}(q) J(q)\dot{q}\leq\sigma_{\max}(\mathcal{Q}(q)J(q))^4\sigma_{\max}(\mathcal{Z}(\dot{\Theta}(t)))^2 \|\dot{q}\|^2$ for all $q,\dot{q}\in\mathbb{R}^n$ from~\eqref{MassMatrix}. Together with the upper bound $\sigma_{\max}(\mathcal{Z}(\dot{\Theta}(t)))\leq \sqrt{69 N}\mu$ from above, this constitutes $\sigma_{\max}(M(q,\dot{\Theta}(t)))\leq \sqrt{69 N}\mu\sigma_{\max}(\mathcal{Q}(q)J(q))^2$, as claimed in~\eqref{BoundedEffectofVaryingParam}.
\end{proof}

Since internal movement of mass-carrying particles triggers time dependency of inertial parameters, realistic robustness testing should also take the parasitic effects of the lumped particle velocity $\Psi$ and of the lumped particle acceleration $\dot{\Psi}$ on the dynamical behavior into account as represented in the generalized robotics equation~\eqref{RoboticEquation} by $H(q,\Psi)\dot{q}$ and $J(q)^\top\mathcal{Q}(q)^\top\dot{\Psi}$, respectively. First, given boundedness of the lumped acceleration with $\|\dot{\Psi}(t)\|\leq \mathrm{d}\Psi^\star$ for all $t\geq 0$, where $\mathrm{d}\Psi^\star\in[0,\infty)$ is some constant, then the upper bound
\begin{align}
\label{LumpedAccelerationBounding}
\|J(q)^\top\mathcal{Q}(q)^\top\dot{\Psi}(t)\|\leq \mathrm{d}\Psi^\star\sup_{q\in\mathbb{R}^n}\big\{\sigma_{\max}(\mathcal{Q}(q)J(q))\big\},
\end{align}
which is valid for all $q\in\mathbb{R}^n$, $t\geq 0$, allows bounding the dynamical effect of internal acceleration of mass-carrying particles under fulfillment of Assumption~\ref{Assumption:UniformUpperBound} as this assumption provides the means to bound the term on the right-hand-side of~\eqref{LumpedAccelerationBounding} with $\sup_{q\in\mathbb{R}^n}\big\{\sigma_{\max}(\mathcal{Q}(q)J(q))\big\}<\infty$, cf.\ Lemma~\ref{lemma:UpperBoundedJacobian}. Therefore, we rely on the same property of the kinematic layout to bound the effects of time-dependent inertial parameters and of acceleration of mass-carrying particles, namely Assumption~\ref{Assumption:UniformUpperBound}---i.e., provided the origins of the body-fixed frames are located at the joints of the robot, we need fulfillment of the specifications listed in~\cite{ghorbel_uniform_1998} for the existence of a finite upper uniform bound of the mass matrix in the case of constant, physically consistent inertial parameters. Next, in order to bound the second parasitic effect $H(q,\Psi)\dot{q}$ on the dynamical behavior that originates from the internal velocities of mass-carrying particles, we need validity of another property of the kinematic layout regarding the partial derivatives of $\mathcal{Q}(\cdot)J(\cdot)$ projected in the specific direction of the lumped velocity $\Psi$ as the matrix $H(q,\Psi)$ represents the skew-symmetric part of ${\partial( \mathcal{Q}(q)J(q))^\top \Psi}/{\partial q}$, cf.~\eqref{DefH}. Let
\begin{align}
\label{DefK}
    \mathcal{K}_p(q)=(\mathcal{Q}(q)J(q))^\top \mathrm{e}_{p,6N}\in\mathbb{R}^{n}, 
\end{align}
$p\in\{1,\ldots,6N\}$, where the standard unit vector $\mathrm{e}_{p,6N}\in\mathbb{R}^{6N}$ picks the $p$th row out of representation $\mathcal{Q}(\cdot)J(\cdot)\in\mathbb{R}^{6N\times n}$ of the kinematic layout of the robotic manipulator, and consider:
\begin{assumption}\label{assumption:BoundedPartialDerivatives}
     $\mathcal{K}_p(\cdot)$ admit $\sup_{q\in\mathbb{R}^n}\Big\{\sigma_{\max}\Big(\tfrac{\partial \mathcal{K}_{p}(q)}{\partial q}-\Big(\tfrac{\partial \mathcal{K}_{p}(q)}{\partial q}\Big)^\top\Big)\Big\}<\infty$ for all $p\in\{1,\ldots,6N\}$.
\end{assumption}
Finally, under validity of Assumption~\ref{assumption:BoundedPartialDerivatives} on the kinematic layout, we can also bound the dynamical effect arising from internal velocity of mass-carrying particles:
\begin{theorem}
    Let Assumption~\ref{assumption:BoundedPartialDerivatives} apply. If the lumped velocity $\Psi$ is bounded with $\|\Psi(t)\|\leq\Psi^\star$ for all $t\geq 0$ with some $\Psi^\star\in[0,\infty)$, then the dynamical effect of $\Psi$ in~\eqref{RoboticEquation} remains bounded since then
    \begin{align}
    \label{BoundedEffectOfLumpedVel}
        \sigma_{\max}(H(q,\Psi(t)))\leq \sqrt{6N}\,\Psi^\star \underbrace{\max_{p\in\{1,\ldots,6N\}}\Big\{\sup_{q\in\mathbb{R}^n}\Big\{\sigma_{\max}\Big(\frac{\partial \mathcal{K}_{p}(q)}{\partial q}-\Big(\frac{\partial \mathcal{K}_{p}(q)}{\partial q}\Big)^\top\Big)\Big\}\Big\}}_{\text{finite as per Assumption~\ref{assumption:BoundedPartialDerivatives}}}
    \end{align}
    for all $q\in\mathbb{R}^n$, $t\geq 0$.
\end{theorem}
\begin{proof}
    Statements in this proof are valid for all $p\in\{1,\ldots,6N\}$, $q\in\mathbb{R}^n$, $t\geq 0$. Consider the decomposition $\Psi(t)=(\begin{smallmatrix}
        \Psi_1(t)&\ldots&\Psi_{6N}(t)
    \end{smallmatrix})^\top$ of the lumped velocity into $\Psi_p(t)\in\mathbb{R}$. Then, from~\eqref{DefH} and from the consequence $J(\cdot)^\top\mathcal{Q}(\cdot)^\top=(\begin{smallmatrix}
        \mathcal{K}_1(\cdot)&\ldots&\mathcal{K}_{6N}(\cdot)
    \end{smallmatrix})$ of~\eqref{DefK}, we obtain $H(q,\Psi(t))=(\sum_{p=1}^{6N}(\partial \mathcal{K}_p(q)/\partial q)\Psi_p(t))-(\sum_{p=1}^{6N}(\partial \mathcal{K}_p(q)/\partial q)\Psi_p(t))^\top=\sum_{p=1}^{6N}Y_p(q)\Psi_p(t)$ with $Y_p(q)= \partial \mathcal{K}_p(q)/\partial q - (\partial \mathcal{K}_p(q)/\partial q)^\top\in\mathbb{R}^{n\times n}$. We get $\dot{q}^\top H(q,\Psi(t))^\top H(q,\Psi(t)) \dot{q}=\|\sum_{p=1}^{6N}Y_p(q)\dot{q}\Psi_p(t)\|^2\leq(\sum_{p=1}^{6N}\|Y_p(q)\dot{q}\| |\Psi_p(t)|)^2\leq (\sum_{p=1}^{6N} \sigma_{\max}(Y_p(q))\cdot\allowbreak|\Psi_p(t)|)^2\|\dot{q}\|^2$ for all $\dot{q}\in\mathbb{R}^n$. This reveals $\sigma_{\max}(H(q,\Psi(t)))\leq \sum_{p=1}^{6N} \sigma_{\max}(Y_p(q))|\Psi_p(t)|$ such that an application of the Cauchy-Schwarz inequality gives $\sigma_{\max}(H(q,\Psi(t)))\leq (\sum_{p=1}^{6N}\sigma_{\max}(Y_p(q))^2)^{1/2}\|\Psi(t)\|\leq \sqrt{6N}\max_{p\in\{1,\ldots,6N\}}\big\{\sigma_{\max}(Y_p(q))\big\}\Psi^\star$. Taking the supremum over all $q\in\mathbb{R}^n$ of the latter upper bound of $\sigma_{\max}(H(q,\Psi(t)))$ results in the desired inequality~\eqref{BoundedEffectOfLumpedVel} and thereby concludes this proof.
\end{proof}

\section{Conclusion and related work}\label{section:Conclusion}
The contribution of this work is threefold. First, we have provided a rigorous proof for the characterization of physical consistency through positive definiteness of a certain symmetric $4\times 4$ matrix-arrangement of the inertial parameters so that other, subsequent findings in this work that utilize this characterization rest on firm, interpretable grounds. Specifically, it turned out that restricting inertial parameters with $m>0$ to originate exclusively via Lebesgue integration from nonnegative mass densities is sufficient to ensure fulfillment of the strict triangle inequalities,\footnote{The proof of sufficiency in the Appendix produces this fact as a byproduct. For details, refer to  Remark~\ref{remark:StrictTriangleIneqsUnderPhysicalConsistency}.} which is crucial for the positive definiteness that characterizes the physical consistency of such parameters. Second, we have derived a generalization of the robotics equation that describes all relevant effects associated with time-dependent change of inertial parameters including these originating from causative time dependencies of the mass densities. The modeling provides a clear distinction between effects that are and those that are not described by this generalization, thereby enabling purposeful testing of e.g.\ adaptive control algorithms by weighing realism and complexity. It is noteworthy that loading processes of end-effectors are inside the scope of our model. Third, we have highlighted kinematic layouts for which uniform physical consistency and upper boundedness of (estimated) inertial parameters guarantees finite, positive uniform boundedness of the (estimated) mass matrix. Now, as is shown in \cite{kaufmann_NICL_2026}---wherein we follow up on the structural findings presented here by leveraging them for a closed-loop analysis---natural adaptation schemes inspired by \cite{lee_natural_2018} in fact already provide uniformity and upper boundedness of their physically consistent estimates. This is a powerful result: For robots with a normal, bounded kinematic layout, it elevates the natural adaptation approach proposed in \cite{lee_natural_2018} from promising in some sense refined control action to providing the guarantee of uniform positive boundedness of the estimated mass matrix, thereby removing a pain point of robot control when only imprecise knowledge about the inertial parameters is available. Therefore, future work will attempt to robustify the natural adaptation framework to the degree that it withstands disturbances as well as the effects of time-varying inertial parameters, even if these parameter variations originate from internal movement of mass-carrying particles as represented herein by the generalized robotics equation.

\section*{Appendix}

\subsection*{Proof of sufficiency in the equivalence claimed in Theorem~\ref{theorem:PhysicalConsistency}}
We want to show that physical consistency of $\Phi\in\mathbb{R}^{10}$ with entries as in~\eqref{eq:Phi} results in $f(\Phi)\succ 0$. To this end, suppose that $\Phi=(\begin{smallmatrix}
        m  & h^\top  & I_{{11}}  & I_{{22}}  & I_{{33}}  & I_{{12}}  & I_{{23}}  & I_{{13}}
    \end{smallmatrix})^{ \top}\in\mathbb{R}^{10}$, i.e.\ structured as per~\eqref{eq:Phi}, is physically consistent. According to our refined Definition~\ref{definition:PhysicalConsistency} of physical consistency, then $m>0$ and we find a function $\rho:\mathbb{R}^3\to[0,\infty)$ that satisfies the equations in~\eqref{MassDistributionToInertialParam} with the integrals therein evaluated with respect to the volumetric Lebesgue measure. We aim our subsequent effort at verifying the strict inequality
    \begin{align}
    \label{GoalForSufficiency}
        \frac{1}{2}\mathrm{tr}(I)&>\lambda_{\max}\left(I+m^{-1}hh^\top\right)
    \end{align}
    since its validity, due to $m>0$, would result in the desired positive definiteness of $f(\Phi)$, see: 
    \begin{align}
        \label{Intermediate}\eqref{GoalForSufficiency}\quad\Longleftrightarrow\quad\frac{1}{2}\mathrm{tr}(I)\mathrm{I}_{3}\succ I+m^{-1}hh^\top \quad\overset{\text{\eqref{DefTrafoInertialParameters}}}{\Longleftrightarrow}\quad\Sigma-m^{-1}hh^\top \succ 0\underset{m>0\text{, \eqref{DefTrafoInertialParameters}}}{\overset{\text{Schur's complement,}}{\Longleftrightarrow}} f(\Phi)\succ 0.
    \end{align}
Exploiting the geometric fact that
\begin{align}
\label{Tjaja}
\mathcal{S}(x)^\top\mathcal{S}(x)=&\mathrm{tr}(xx^\top)\mathrm{I}_3-xx^\top
\end{align}
holds for all $x\in\mathbb{R}^3$, we rewrite the term on the left-hand-side of~\eqref{GoalForSufficiency} as
\begin{align}
    \frac{1}{2}\mathrm{tr}(I)\overset{\eqref{MassDistributionToInertialParam}}{=}&\frac{1}{2}\int_{\mathbb{R}^3}\rho(x)\mathrm{tr}\left(\mathcal{S}(x)^\top\mathcal{S}(x)\right)\mathrm{d}x\overset{\eqref{Tjaja}}{=}\frac{1}{2}\int_{\mathbb{R}^3}\rho(x)\mathrm{tr}\left(\mathrm{tr}\left(xx^\top\right)\mathrm{I}_3-xx^\top\right)\mathrm{d}x\allowdisplaybreaks\\\label{LHS}
    =&\int_{\mathbb{R}^3}\rho(x)\mathrm{tr}\left(xx^\top\right)\mathrm{d}x
\end{align}
and the term on its right-hand-side as
\begin{align}\label{Intermediate1}
    \lambda_{\max}\left(I+m^{-1}hh^\top\right)=&\max_{v\in\mathbb{R}^3,\ \|v\|=1}\left\{v^\top \left(I+m^{-1}hh^\top\right)v\right\}\allowdisplaybreaks\\\overset{\eqref{MassDistributionToInertialParam}}{=}&\max_{v\in\mathbb{R}^3,\ \|v\|=1}\left\{\int_{\mathbb{R}^3}\rho(x)v^\top \mathcal{S}(x)^\top\mathcal{S}(x)v\mathrm{d}x+m^{-1}\left(v^\top h\right)^2\right\}\allowdisplaybreaks\\
    \overset{\eqref{Tjaja}}{=}&\max_{v\in\mathbb{R}^3,\ \|v\|=1}\left\{\int_{\mathbb{R}^3}\rho(x)v^\top \left(\mathrm{tr}\left(xx^\top\right)\mathrm{I}_3-xx^\top\right)v\mathrm{d}x+m^{-1}\left(v^\top h\right)^2\right\}\allowdisplaybreaks\\
    \label{Intermediate2}=&\int_{\mathbb{R}^3}\rho(x)\mathrm{tr}\left(xx^\top\right)\mathrm{d}x+\max_{v\in\mathbb{R}^3,\ \|v\|=1}\left\{m^{-1}\left(v^\top h\right)^2-\int_{\mathbb{R}^3}\rho(x) \left(v^\top x\right)^2\mathrm{d}x \right\}.
\end{align}
Thereby, we arrive at
\begin{align}
    \label{Intermediate0}\eqref{GoalForSufficiency}\quad\Longleftrightarrow\quad& 0>\max_{v\in\mathbb{R}^3,\ \|v\|=1}\left\{m^{-1}\left(v^\top h\right)^2-\int_{\mathbb{R}^3}\rho(x) \left(v^\top x\right)^2\mathrm{d}x \right\}\allowdisplaybreaks\\
    \quad\overset{m>0}{\Longleftrightarrow}\quad& \forall v\in\mathbb{R}^3\text{ with }\|v\|=1: \left(v^\top h\right)^2 <m\int_{\mathbb{R}^3}\rho(x) \left(v^\top x\right)^2\mathrm{d}x\allowdisplaybreaks\\
    \quad\overset{\eqref{MassDistributionToInertialParam}}{\Longleftrightarrow}\quad& \forall v\in\mathbb{R}^3\text{ with }\|v\|=1:\left(\int_{\mathbb{R}^3}\rho(x)v^\top x\mathrm{d}x\right)^2<\int_{\mathbb{R}^3}\rho(x)\mathrm{d}x\int_{\mathbb{R}^3}\rho(x) \left(v^\top x\right)^2\mathrm{d}x\allowdisplaybreaks\\
    \quad{\Longleftrightarrow}\quad& \forall v\in\mathbb{R}^3\text{ with }\|v\|=1: \left(\int_{\mathbb{R}^3}\mathscr{f}(x)\mathscr{g}_v(x)\mathrm{d}x\right)^2<\int_{\mathbb{R}^3}\mathscr{f}(x)^2\mathrm{d}x\int_{\mathbb{R}^3}\mathscr{g}_v(x)^2\mathrm{d}x,\label{StrictSchwartz}
\end{align}
where $\mathscr{f}(x)=\sqrt{\rho(x)}$ and $\mathscr{g}_v(x)=\sqrt{\rho(x)}v^\top x$ for $x\in\mathbb{R}^3$. In order to reformulate the statement on the right-hand-side of~\eqref{StrictSchwartz} by means of the strict Schwarz inequality, we need to ensure its applicability by confirming $\mathscr{f},\mathscr{g}_v\in{L}^2(\mathbb{R}^3)$ for all $v\in\mathbb{R}^3$ with $\|v\|=1$. First, $\mathscr{f}\in{L}^2(\mathbb{R}^3)$ holds due to $\|\mathscr{f}\|_{{L}^2}=\sqrt{m}<\infty$. Second, for all $v\in\mathbb{R}^3$ with $\|v\|=1$, we obtain
\begin{align}
    \|\mathscr{g}_v\|_{{L}^2}^2=&\int_{\mathbb{R}^3}\rho(x) v^\top xx^\top v\mathrm{d}x\overset{\eqref{Tjaja}}{=}\int_{\mathbb{R}^3}\rho(x)  \left(\mathrm{tr}\left(xx^\top\right)-v^\top\mathcal{S}(x)^\top\mathcal{S}(x)v\right)\mathrm{d}x\overset{\eqref{LHS},\ \eqref{MassDistributionToInertialParam}}{=}\frac{1}{2}\mathrm{tr}(I)-v^\top I v\allowdisplaybreaks\\
    \in&\left[\frac{1}{2}\mathrm{tr}(I)-\lambda_{\max}(I),\ \frac{1}{2}\mathrm{tr}(I)-\lambda_{\min}(I)\right]\overset{\eqref{DefTrafoInertialParameters}}{=}\left[\lambda_{\min}(\Sigma),\ \lambda_{\max}(\Sigma)\right],\label{Intermediate3}
\end{align}
wherein the lower bound for $\|\mathscr{g}_v\|_{{L}^2}^2$ with $\lambda_{\min}(\Sigma)\geq 0$ is nonnegative since we have
\begin{align}
    \Sigma\overset{\eqref{DefTrafoInertialParameters}}{=}&\frac{1}{2}\mathrm{tr}(I)\mathrm{I}_3-I\overset{\eqref{LHS},\ \eqref{MassDistributionToInertialParam}}{=}\int_{\mathbb{R}^3}\rho(x)\mathrm{tr}\left(xx^\top\right)\mathrm{I}_3-\rho(x)\mathcal{S}(x)^\top\mathcal{S}(x)\mathrm{d}x\overset{\eqref{Tjaja}}{=}\int_{\mathbb{R}^3}\rho(x)xx^\top\mathrm{d}x
\end{align}
such that $\Sigma$ inherits positive semi-definiteness from $\rho(x)xx^\top\succeq 0$ which is given for all $x\in\mathbb{R}^3$ as $\rho$ is nonnegative. Thereby,~\eqref{Intermediate3} results in the desired $\mathscr{g}_v\in{L}^2(\mathbb{R}^3)$ for all $v\in\mathbb{R}^3$ with $\|v\|=1$ since $\|\mathscr{g}_v\|_{{L}^2}\in\mathcal{I}$ with the nonnegative interval $\mathcal{I}=\big[\sqrt{\lambda_{\min}(\Sigma)},\ \sqrt{\lambda_{\max}(\Sigma)}\big]$ that is bounded above since $\Phi\in\mathbb{R}^{10}$ yields $\Sigma\in\mathbb{R}^{3\times 3}$, cf.~\eqref{DefTrafoInertialParameters}. Further, in addition to their square-integrability, the nonnegativity of $\rho$ also guarantees that the functions $\mathscr{f},\mathscr{g}_v$ for all $v\in\mathbb{R}^3$ with $\|v\|=1$ map exclusively into the real numbers, thus revealing them as admissible arguments for $\langle\cdot,\cdot\rangle:({L}^2(\mathbb{R}^3))^2\to\mathbb{R}$ with $\langle\mathscr{f},\mathscr{g}_v\rangle=\int_{\mathbb{R}^3}\mathscr{f}(x)\mathscr{g}_v(x)\mathrm{d}x$ which is a scalar product for the vector space of square-integrable, real-valued functions defined in $3$-dimensional space. This renders the strict Schwarz inequality (see e.g.\  (\cite{folland_real_1999}, $5.19$)) applicable to find an upper bound of $|\langle\mathscr{f},\mathscr{g}_v\rangle|$ and therefore,~\eqref{Intermediate0}-\eqref{StrictSchwartz} yields
\begin{align}
    \label{Intermediate5}\eqref{GoalForSufficiency}\quad\Longleftrightarrow\quad&\forall v\in\mathbb{R}^3\text{ with }\|v\|=1:|\langle\mathscr{f},\mathscr{g}_v\rangle|<\sqrt{\langle\mathscr{f},\mathscr{f}\rangle}\ \sqrt{\langle\mathscr{g}_v,\mathscr{g}_v\rangle}\allowdisplaybreaks\\
    \overset{\text{(\cite{folland_real_1999}, $5.19$)}}{\Longleftrightarrow}\quad&\forall v\in\mathbb{R}^3\text{ with }\|v\|=1: \neg\underbrace{\left(\exists 0\neq\begin{pmatrix}\delta_1& \delta_2\end{pmatrix}^\top\in\mathbb{R}^2: \left\langle\delta_1\mathscr{f}+\delta_2\mathscr{g}_v,\delta_1\mathscr{f}+\delta_2\mathscr{g}_v\right\rangle=0\right)}_{\text{linear dependence of $\mathscr{f}$, $\mathscr{g}_v$ with respect to\ $\langle\cdot,\cdot\rangle$}}\allowdisplaybreaks\\
    {\Longleftrightarrow}\quad&\forall v\in\mathbb{R}^3\text{ with }\|v\|=1: \neg\left(\exists 0\neq\begin{pmatrix}\delta_1& \delta_2\end{pmatrix}^\top\in\mathbb{R}^2: \delta_1\mathscr{f}+\delta_2\mathscr{g}_v=0\text{ a.e.}\right)\allowdisplaybreaks\\
    {\Longleftrightarrow}\quad&\forall v\in\mathbb{R}^3\text{ with }\|v\|=1: \neg\left(\exists 0\neq\begin{pmatrix}\delta_1& \delta_2\end{pmatrix}^\top\in\mathbb{R}^2: \mu\left(\mathscr{M}_v(\delta_1,\delta_2)\right)=0\right)\allowdisplaybreaks\\
    {\Longleftrightarrow}\quad&\forall v\in\mathbb{R}^3\text{ with }\|v\|=1,\   0\neq\begin{pmatrix}\delta_1& \delta_2\end{pmatrix}^\top\in\mathbb{R}^2: \mu\left(\mathscr{M}_v(\delta_1,\delta_2)\right)>0.\label{Intermediate4}
\end{align}
Therein, we applied {the Lebesgue measure $\mu$ for volumina in $\mathbb{R}^3$} to the set $\mathscr{M}_v(\delta_1,\delta_2)=\big\{x\in\mathbb{R}^3:\delta_1\mathscr{f}(x)+\delta_2\mathscr{g}_v(x)\neq 0\big\}$. Note that, due to the fact that $\mathscr{f}(x),\mathscr{g}_v(x)= 0$ if and only if $x\notin\mathrm{carr}(\rho)$, where $\mathrm{carr}(\rho)=\big\{x\in\mathbb{R}^3:\rho(x)\neq 0\big\}$ denotes the \textit{carrier} of $\rho$, we can rewrite $\mathscr{M}_v(\delta_1,\delta_2)$ as
\begin{align}
    \!\!\!\!\mathscr{M}_v(\delta_1,\delta_2)=&\left\{x\in\mathrm{carr}(\rho):\delta_1\mathscr{f}(x)+\delta_2\mathscr{g}_v(x)\neq 0\right\}=\left\{x\in\mathrm{carr}(\rho):\rho(x)\begin{pmatrix}1&v^\top x\end{pmatrix}\begin{pmatrix}
        \delta_1\\\delta_2
    \end{pmatrix}\neq 0\right\}
    \allowdisplaybreaks\\
    =&\left\{x\in\mathrm{carr}(\rho):\begin{pmatrix}1&v^\top x\end{pmatrix}\begin{pmatrix}
        \delta_1\\\delta_2
    \end{pmatrix}\neq 0\right\}.\label{SetRepresentation}
\end{align}
Using the representation of $\mathscr{M}_v(\delta_1,\delta_2)$ in~\eqref{SetRepresentation}, we are now in a position to validate the strict inequality \eqref{GoalForSufficiency} and hence, to obtain the desired $f(\Phi)\succ 0$ by proving the statement on the right hand-side of~\eqref{Intermediate4}, cf.~\eqref{Intermediate4}-\eqref{Intermediate5},~\eqref{Intermediate}. To this end, we take any $v\in\mathbb{R}^3$ with $\|v\|=1$,  $  0\neq\Big(\begin{smallmatrix}\delta_1\\ \delta_2\end{smallmatrix}\Big)\in\mathbb{R}^2$ and aim to infer that then $\mu\left(\mathscr{M}_v(\delta_1,\delta_2)\right)>0$. 
Note that $m>0$ together with the nonnegativity of $\rho$, by virtue of (\cite{folland_real_1999}, Proposition $2.16$), prohibits the equality $\rho=0$ to be satisfied a.e., or, equivalently, demands $\mu(\mathrm{carr}(\rho))>0$. This fact is useful to assess the measure of $\mathscr{M}_v(\delta_1,\delta_2)$ in either of the subsequent scenarios, where one of them certainly occurs:
\begin{enumerate}
    \item $\mathscr{M}_v(\delta_1,\delta_2)=\mathrm{carr}(\rho)$. Then, we have the desired $f(\Phi)\succ 0$ as $\mu\left(\mathscr{M}_v(\delta_1,\delta_2)\right)>0$ is immediate.
    \item $\mathscr{M}_v(\delta_1,\delta_2)\subset\mathrm{carr}(\rho)$. That is, the carrier of $\rho$ is not fully occupied by $\mathscr{M}_v(\delta_1,\delta_2)$ and hence, we find $x_1\in \mathrm{carr}(\rho)\setminus \mathscr{M}_v(\delta_1,\delta_2)$ since $\mathrm{carr}(\rho)$ is not empty due to its positive measure. Consider the hyperplane $\mathscr{H}_v=\{x\in\mathbb{R}^3:v^\top x=0\}$ perpendicular to $v$. Since this hyperplane, shifted into $\mathrm{carr}(\rho)$ and placed at $x_1$, is a null set, i.e.\ $\mu(\mathscr{H}_v\oplus\{x_1\})=0$, we obtain $\mu(\mathscr{N}_v(x_1))>0$ for the set $\mathscr{N}_v(x_1)=\mathrm{carr}(\rho)\setminus(\mathscr{H}_v\oplus \{x_1\})$ denoting the carrier of $\rho$ without the shifted hyperplane. According to its positive measure, $\mathscr{N}_v(x_1)$ is not empty and, due to its construction, any $x_2\in\mathscr{N}_v(x_1)$ admits $x_2\in\mathrm{carr}(\rho)$ as well as $x_2-x_1\notin \mathscr{H}_v$, i.e.\ $v^\top(x_2-x_1)\neq 0$, where the latter also reads as $\mathrm{det}\Big(\begin{smallmatrix}
    1&v^\top x_2\\1&v^\top x_1
\end{smallmatrix}\Big)\neq 0$ such that we obtain $\Big(\begin{smallmatrix}
    1&v^\top x_2\\1&v^\top x_1
\end{smallmatrix}\Big)\Big(\begin{smallmatrix}\delta_1\\\delta_2\end{smallmatrix}\Big)\neq 0$ from $\Big(\begin{smallmatrix}\delta_1\\\delta_2\end{smallmatrix}\Big)\neq 0$. Now, $(\begin{smallmatrix}
    1&v^\top x_1
\end{smallmatrix})\Big(\begin{smallmatrix}\delta_1\\\delta_2\end{smallmatrix}\Big)=0$ is a consequence of $x_1\notin M_v(\delta_1,\delta_2)$,~cf.\ \eqref{SetRepresentation}. Thereby, we are left with $(\begin{smallmatrix}
    1&v^\top x_2
\end{smallmatrix})\Big(\begin{smallmatrix}\delta_1\\\delta_2\end{smallmatrix}\Big)\neq 0$ such that $x_2\in\mathrm{carr}(\rho)$, in view of~\eqref{SetRepresentation}, necessitates $x_2\in\mathscr{M}_v(\delta_1,\delta_2)$. Hence, we have $\mathscr{N}_v(x_1)\subseteq \mathscr{M}_v(\delta_1,\delta_2)$ and therefore obtain positivity of the measure of $\mathscr{M}_v(\delta_1,\delta_2)$ via $\mu(\mathscr{M}_v(\delta_1,\delta_2))\geq\mu(\mathscr{N}_v(x_1))>0$ by means of the monotonicity of $\mu$, thus also leading to the desired $f(\Phi)\succ 0$.
\end{enumerate}
This concludes the proof of the sufficient direction in the equivalence claimed in Theorem~\ref{theorem:PhysicalConsistency}.
\begin{remark}\label{remark:StrictTriangleIneqsUnderPhysicalConsistency}
    Above, we have shown that physical consistency of inertial parameters $\Phi\in\mathbb{R}^{10}$ in accordance with Definition~\ref{definition:PhysicalConsistency} ensures validity of the strict inequality~\eqref{GoalForSufficiency}. While the latter, due to~\eqref{Intermediate} and the required $m>0$, is equivalent to the desired $f(\Phi)\succ 0$, it also guarantees $\Sigma_\mathrm{CoM}=\Sigma-m^{-1}hh^\top\succ 0$, which as per~\eqref{SigmaCoMPosDef_IFF_StricTriangleCoM} guarantees the inertia matrix $I_\mathrm{CoM}=I-m^{-1}\mathcal{S}(h)^\top\mathcal{S}(h)$ for rotation about the CoM of the rigid body to fulfill the strict triangle inequalities. That is, even though not explicitly mentioned in our proof of the sufficient direction above, it closes the gap that was left by~\cite{wensing_linear_2018}: Ambiguity about whether or not these triangle inequalities are strict under physical consistency of the inertial parameters.
\end{remark}

\subsection*{Derivation of the generalized robotics equation by using the Lagrange formalism}\label{Appendix:RoboticEquation}
In order to shorten some expressions in the subsequent analysis, we use the following convention: Whenever $J$, $M$, $M_{k,j}$, $\mathcal{Q}$, $U$, $L$,  $\Gamma_{i,j,k}$,  $C$, $\mathcal{F}$, $G$ are written without arguments, they represent $J(q)$, $M(q,\Theta)$, $M_{k,j}(q,\Theta)$, $\mathcal{Q}(q)$, $U(q,\Theta)$, $L(q,\dot{q},\Theta,t)$,  $\Gamma_{i,j,k}(q,\Theta)$,  $C(q,\dot{q},\Theta)$, $\mathcal{F}(q,\Theta,\dot{\Theta})$, $G(q,\Theta)$, respectively.

We use~\eqref{KineticEnergy},~\eqref{Lagrangian} and~\eqref{LagrangeFormalism} to derive the generalization~\eqref{RoboticEquation} of the robotics equation. To this end, let $k\in\{1,\ldots,n\}$ and consider the decomposition of $J=(\begin{smallmatrix}
    J_1&\ldots&J_n
\end{smallmatrix})\in\mathbb{R}^{6N\times n}$ into $J_k\in\mathbb{R}^{6N}$. Due to $M\in\mathrm{sym}(n)$ and $\dot{q}^\top J^\top \mathcal{Q}^\top=\sum_{j=1}^n \dot{q}_j J_j^\top\mathcal{Q}^\top $, we get 
\begin{align}
\frac{\partial L}{\partial \dot{q}_k}&=\sum_{j=1}^n \Big(\frac{1}{2}\Big(\sum_{i=1}^n M_{i,j}\frac{\partial \dot{q}_i\dot{q}_j}{\partial \dot{q}_k}\Big)+\frac{\partial  \dot{q}_jJ_j^\top\mathcal{Q}^\top \Psi}{\partial \dot{q}_k}\Big)\allowdisplaybreaks\\
&=\Big(\sum_{j=1}^n M_{k,j}\dot{q}_j\Big)+J_k^\top\mathcal{Q}^\top\Psi
\end{align}
and thus obtain
\begin{align}
    \frac{\mathrm{d}}{\mathrm{d}t}\Big(\frac{\partial L}{\partial \dot{q}_k}\Big)&=\Big(\sum_{j=1}^n\frac{\mathrm{d} M_{k,j}\dot{q}_j}{\mathrm{d}t}\Big)+\frac{\mathrm{d}\Psi^\top \mathcal{Q} J_k}{\mathrm{d}t}\allowdisplaybreaks\\
    &=\Big(\sum_{j=1}^n M_{k,j}\ddot{q}_j+\Big(\sum_{i=1}^n\frac{\partial M_{k,j}}{\partial q_i}\dot{q}_i\dot{q}_j\Big)+\Big(\sum_{h=1}^{10N}\frac{\partial M_{k,j}}{\partial \Theta_h}\dot{\Theta}_h\dot{q}_j\Big)\Big)+\dot{\Psi}^\top \mathcal{Q} J_k+\frac{\partial\Psi^\top \mathcal{Q} J_k}{\partial q}\dot{q}
\end{align}
from $J_k^\top\mathcal{Q}^\top \Psi=\Psi^\top \mathcal{Q} J_k$. Further, the structure of the Lagrangian in~\eqref{Lagrangian} leads to
\begin{align}
    \frac{\partial L}{\partial q_k}&=\Big(\sum_{j=1}^n\Big(\frac{1}{2}\sum_{i=1}^n\frac{\partial M_{i,j}}{\partial q_k}\dot{q}_i\dot{q}_j\Big)+\frac{\partial \dot{q}_j J_j^\top\mathcal{Q}^\top \Psi}{\partial q_k}\Big)-\frac{\partial U}{\partial q_k}
\end{align}
such that the equations of motion~\eqref{LagrangeFormalism} of the open kinematic chain evaluate to
\begin{align}
\label{ExplicitLagrange}
    &\Big(\sum_{j=1}^n M_{k,j}\ddot{q}_j\Big)+\Gamma +\Big(\sum_{j=1}^n\sum_{h=1}^{10N} \frac{\partial M_{k,j}}{\partial \Theta_h}\dot{\Theta}_h\dot{q}_j\Big)+\frac{\partial U}{\partial q_k}\\&-\Big(\sum_{j=1}^n\frac{\partial   \Psi^\top\mathcal{Q} J_j}{\partial q_k}\dot{q}_j\Big)+J_k^\top\mathcal{Q}^\top \dot{\Psi} +\frac{\partial J_k^\top\mathcal{Q}^\top \Psi}{\partial q}\dot{q}=\tau_k+w_k\nonumber,
\end{align}
where
\begin{align}
    \Gamma=&\sum_{j=1}^n\sum_{i=1}^n\left(\frac{\partial M_{k,j}}{\partial q_i}-\frac{1}{2}\frac{\partial M_{i,j}}{\partial q_k}\right)\dot{q}_i\dot{q}_j=\sum_{j=1}^n\sum_{i=1}^n \Gamma_{i,j,k}\dot{q}_i\dot{q}_j
\end{align}
in view of the definition~\eqref{Cristoffel} of the Christoffel symbols of the first kind. Thus, the equation~\eqref{ExplicitLagrange} can be rewritten as
\begin{align}\label{skalarRoboticEquation}
    &\Big(\sum_{j=1}^n M_{k,j}\ddot{q}_j\Big)+\Big(\sum_{j=1}^n\Big(\sum_{i=1}^n\Gamma_{i,j,k}\dot{q}_i+\sum_{h=1}^{10N} \frac{\partial M_{k,j}}{\partial \Theta_h}\dot{\Theta}_h\Big)\dot{q}_j\Big)\\&
    +\Big(\frac{\partial J_k^\top\mathcal{Q}^\top \Psi}{\partial q}-\frac{\partial \Psi^\top\mathcal{Q} J }{\partial q_k}\Big)\dot{q}+\frac{\partial U}{\partial q_k}= \tau_k+w_k-J_k^\top\mathcal{Q}^\top \dot{\Psi}\nonumber.
\end{align}
Stacking terms on top of each other as in
\begin{align}
    \begin{pmatrix}
        \frac{\partial J_1^\top \mathcal{Q}^\top \Psi}{\partial q}-\frac{\partial \Psi^\top \mathcal{Q}  J }{\partial q_1}\\
        \vdots\\
        \frac{\partial J_n^\top \mathcal{Q}^\top \Psi}{\partial q}-\frac{\partial \Psi^\top \mathcal{Q}  J }{\partial q_n}
    \end{pmatrix}&=\frac{\partial J^\top\mathcal{Q}^\top\Psi}{\partial q}-\begin{pmatrix}
        \Big(\frac{\partial J^\top \mathcal{Q}^\top\Psi}{\partial q_1}\Big)^{ \top}\\
        \vdots\\
        \Big(\frac{\partial J^\top \mathcal{Q}^\top\Psi}{\partial q_n}\Big)^{ \top}
    \end{pmatrix}\allowdisplaybreaks\\
    &=\frac{\partial J^\top \mathcal{Q}^\top\Psi}{\partial q}-\Big(\frac{\partial J^\top \mathcal{Q}^\top\Psi}{\partial q}\Big)^{ \top}
\end{align}
reveals the stacked version of~\eqref{skalarRoboticEquation} for $k=1,\ldots,n$ as 
\begin{align}
\label{IntermediateRoboticEquation}
    M(q,\Theta)\ddot{q} + \left(C(q,\dot{q},\Theta)+\mathcal{F}(q,\Theta,\dot{\Theta})+H(q,\Psi)\right)\dot{q}+G(q,\Theta)=\tau+w-J(q)^\top \mathcal{Q}(q)^\top\dot{\Psi},
\end{align}
when taking into account the definitions~\eqref{DefCoriolis},~\eqref{DefGravitation} and~\eqref{DefH}. In accordance with~\eqref{skalarRoboticEquation}, the matrix $\mathcal{F}(q,\Theta,\dot{\Theta})\in\mathrm{sym}(n)$ in~\eqref{IntermediateRoboticEquation} reads
\begin{align}
\label{DefF}
\mathcal{F}(q,\Theta,\dot{\Theta})&=\left(\sum_{h=1}^{10N}\frac{\partial M_{k,j}(q,\Theta)}{\partial \Theta_h}\dot{\Theta}_h\right)_{(k,j)\in\{1,\ldots,n\}^2},
\end{align}
which can be simplified further, as is discussed next: Since the mass matrix $M(q,\Theta)$ is linear in the inertial parameters, we find regressor functions $M_{k,j}^h(q)\in\mathbb{R}$, $h=\{1,\ldots,10 N\}$ with 
\begin{align}
\label{MassMatrixRegressors}
M_{k,j}(q,\Theta)=\sum_{h=1}^{10N}M_{k,j}^h(q)\Theta_h
\end{align}
for all $k,j\in\{1,\ldots,n\}$. They allow us to rewrite
\begin{align}
    \mathcal{F}(q,\Theta,\dot{\Theta})&=\left(\sum_{h=1}^{10N}\frac{\partial \sum_{\mathscr{h}=1}^{10N}M_{k,j}^\mathscr{h}(q)\Theta_\mathscr{h}}{\partial \Theta_h}\dot{\Theta}_h\right)_{(k,j)\in\{1,\ldots,n\}^2}\allowdisplaybreaks\\
    &=\left(\sum_{h=1}^{10N} M_{k,j}^h(q)\dot{\Theta}_h\right)_{(k,j)\in\{1,\ldots,n\}^2}\allowdisplaybreaks\\
    &=M(q,\dot{\Theta}),
\end{align}
thus shaping the already established~\eqref{IntermediateRoboticEquation} into the desired generalized robotics equation~\eqref{RoboticEquation}.

\subsection*{Verification of two inherent properties of the generalized robotics equation}
We start by verifying the~\ref{Prop1}. inherent property, i.e.\ skew-symmetry of the matrices $\mathcal{N}=\dot{M}-2C-\mathcal{F}$ and $H$. First, $H\in\mathrm{skew}(n)$ is an immediate consequence of its definition in~\eqref{DefH}. Second, due to $\dot{M},\mathcal{F}\in\mathrm{sym}(n)$, we have $\mathcal{N}\in\mathrm{skew}(n)$ if and only if $\mathcal{N}=-\dot{M}+2C^\top+\mathcal{F}$, holding, if and only if $\dot{M}-\mathcal{F}=C+C^\top$. Accordingly, we evaluate the validity of the equation $\dot{M}-\mathcal{F}=C+C^\top$: Due to~\eqref{Cristoffel} and $M\in\mathrm{sym}(n)$, its right-hand-side
\begin{align}
    C+C^\top&=\left(\sum_{i=1}^n\left(\Gamma_{i,j,k}+\Gamma_{i,k,j}\right)\dot{q}_i\right)_{(k,j)\in\{1,\ldots,n\}^2}\allowdisplaybreaks\\
    &=\left(\sum_{i=1}^n\frac{1}{2}\left(\frac{\partial M_{k,j}}{\partial q_i}+\frac{\partial M_{k,i}}{\partial q_j}-\frac{\partial M_{i,j}}{\partial q_k}+\frac{\partial M_{j,k}}{\partial q_i}+\frac{\partial M_{j,i}}{\partial q_k}-\frac{\partial M_{i,k}}{\partial q_j}\right)\dot{q}_i\right)_{(k,j)\in\{1,\ldots,n\}^2}\allowdisplaybreaks\\
    &=\left(\sum_{i=1}^n\frac{\partial M_{k,j}}{\partial q_i}\dot{q}_i\right)_{(k,j)\in\{1,\ldots,n\}^2} 
\end{align}
equals its left-hand-side
\begin{align}
    \dot{M}-\mathcal{F}&=\left(\dot{M}_{k,j}-\Big(\sum_{h=1}^{10N}\frac{\partial M_{k,j}}{\partial \Theta_h}\dot{\Theta}_h\Big)\right)_{(k,j)\in\{1,\ldots,n\}^2}\allowdisplaybreaks\\
    &=\left(\Big(\sum_{i=1}^n \frac{\partial M_{k,j}}{\partial q_i}\dot{q}_i\Big)+\Big(\sum_{h=1}^{10N}\frac{\partial M_{k,j}(q,\Theta)}{\partial \Theta_h}\dot{\Theta}_h\Big)-\Big(\sum_{h=1}^{10N}\frac{\partial M_{k,j}(q,\Theta)}{\partial \Theta_h}\dot{\Theta}_h\Big)\right)_{(k,j)\in\{1,\ldots,n\}^2}.
\end{align}
This results in $\mathcal{N}\in\mathrm{skew}(n)$, as desired.

Next, we back up the~\ref{Prop2}. inherent property of the generalized robotics equation~\eqref{RoboticEquation}, i.e.\ its linear dependence on inertial parameters, on their time derivatives as stated in~\eqref{RoboticRegressor},~\eqref{RoboticRegressorVel}, respectively. Since the mass matrix $M(q,\Theta)$ as well as the potential energy $\mathcal{U}(q,\Theta)$ are linear in the inertial parameters, we find regressor functions $\mathcal{U}^h(q)\in\mathbb{R}$, $h=\{1,\ldots,10 N\}$ with
\begin{align}
\mathcal{U}(q,\Theta)=\sum_{h=1}^{10N} \mathcal{U}^h(q)\Theta_h,
\end{align}
in addition to the already established regressor functions $M_{k,j}^h(q)$ for $h=\{1,\ldots,10 N\}$, $k,j\in\{1,\ldots,n\}$ that constitute the decomposition~\eqref{MassMatrixRegressors}. Accordingly, for any $a=\big(a_k\big)_{k\in\{1,\ldots,n\}},v=\big(v_k\big)_{k\in\{1,\ldots,n\}}\in\mathbb{R}^n$, we obtain $M(q,\Theta)a=R_\mathrm{M}(q,a)\Theta$, $C(q,\dot{q},\Theta)v=R_\mathrm{C}(q,\dot{q},v)\Theta$ and $G(q,\Theta)=R_\mathrm{G}(q)\Theta$ by introducing the regressor functions
\begin{align}
    R_\mathrm{M}(q,a)=&\left(\sum_{j=1}^n M_{k,j}^h(q)a_j\right)_{(k,h)\in\{1,\ldots,n\}\times \{1,\ldots,10N\}},\allowdisplaybreaks\\
    R_\mathrm{C}(q,\dot{q},v)=&\left(\frac{1}{2}\sum_{j=1}^n\sum_{i=1}^n\left(\frac{\partial M_{k,j}^h(q)}{\partial q_i}+\frac{\partial M_{k,i}^h(q)}{\partial q_j}-\frac{\partial M_{i,j}^h(q)}{\partial q_k}\right)\dot{q}_i v_j\right)_{(k,h)\in\{1,\ldots,n\}\times \{1,\ldots,10N\}},\allowdisplaybreaks\\
    R_\mathrm{G}(q)=&\begin{pmatrix}
        \mathcal{U}^1(q)&\ldots&\mathcal{U}^{10N}(q)
    \end{pmatrix}.
\end{align}
Hence, we arrive at  $M (q,\Theta)a + C (q,\dot{q},\Theta)v + G (q,\Theta)
 = R(q,\dot{q},v,a)\Theta$ with $R(q,\dot{q},v,a)=R_\mathrm{M}(q,a)+R_\mathrm{C}(q,\dot{q},v)+R_\mathrm{G}(q)\in\mathbb{R}^{n\times 10N}$, i.e.\ the desired equation~\eqref{RoboticRegressor} results from decomposing $R(q,\dot{q},v,a)=(\begin{smallmatrix}
        R_1(q,\dot{q},v,a)&\ldots&R_{10N}(q,\dot{q},v,a)\end{smallmatrix})$ into $R_l(q,\dot{q},v,a)\in\mathbb{R}^{n\times 10}$, $l\in\{1,\ldots,N\}$. Furthermore, decomposing $R_\mathrm{M}(q,v)=(\begin{smallmatrix}
        R_{\mathrm{M},1}(q,v)&\ldots&R_{\mathrm{M},N}(q,v)\end{smallmatrix})$ into $R_{\mathrm{M},l}(q,v)\in\mathbb{R}^{n\times 10}$, $l\in\{1,\ldots,N\}$ yields~\eqref{RoboticRegressorVel}, as claimed.

\bibliographystyle{ieeetr} 
\bibliography{Bibliothek}

@inproceedings{lee_natural_2018,
	title = {A {Natural} {Adaptive} {Control} {Law} for {Robot} {Manipulators}},
	doi = {10.1109/IROS.2018.8593727},
	booktitle = {{IEEE}/{RSJ} {International} {Conference} on {Intelligent} {Robots} and {Systems}},
	author = {Lee, Taeyoon and Kwon, Jaewoon and Park, Frank C.},
	month = oct,
	year = {2018},
	pages = {1--9},
}

@article{wensing_linear_2018,
	title = {Linear {Matrix} {Inequalities} for {Physically} {Consistent} {Inertial} {Parameter} {Identification}: {A} {Statistical} {Perspective} on the {Mass} {Distribution}},
	volume = {3},
	shorttitle = {Linear {Matrix} {Inequalities} for {Physically} {Consistent} {Inertial} {Parameter} {Identification}},
	doi = {10.1109/LRA.2017.2729659},
	number = {1},
	journal = {IEEE Robotics and Automation Letters},
	author = {Wensing, Patrick M. and Kim, Sangbae and Slotine, Jean-Jacques E.},
	month = jan,
	year = {2018},
	pages = {60--67},
}

@article{slotine_adaptive_1987,
	title = {On the {Adaptive} {Control} of {Robot} {Manipulators}},
	volume = {6},
	issn = {0278-3649},
	doi = {10.1177/027836498700600303},
	language = {en},
	number = {3},
	urldate = {2024-07-17},
	journal = {The International Journal of Robotics Research},
	author = {Slotine, Jean-Jacques E. and Li, Weiping},
	month = sep,
	year = {1987},
	pages = {49--59},
}

@article{craig_adaptive_1987,
	title = {Adaptive {Control} of {Mechanical} {Manipulators}},
	volume = {6},
	issn = {0278-3649},
	doi = {10.1177/027836498700600202},
	language = {en},
	number = {2},
	urldate = {2024-07-17},
	journal = {The International Journal of Robotics Research},
	author = {Craig, John J. and Hsu, Ping and Sastry, S. Shankar},
	month = jun,
	year = {1987},
	pages = {16--28},
}

@article{patnaik_adaptive_2023,
	title = {Adaptive {Attitude} {Control} for {Foldable} {Quadrotors}},
	volume = {7},
	issn = {2475-1456},
	doi = {10.1109/LCSYS.2023.3234045},
	urldate = {2024-07-17},
	journal = {IEEE Control Systems Letters},
	author = {Patnaik, Karishma and Zhang, Wenlong},
	year = {2023},
	pages = {1291--1296},
}

@article{wu_adaptive_2022,
	title = {Adaptive {Tracking} {Control} {With} {Uncertainty}-{Aware} and {State}-{Dependent} {Feedback} {Action} {Blending} for {Robot} {Manipulators}},
	volume = {7},
	issn = {2377-3766},
	doi = {10.1109/LRA.2022.3212669},
	number = {4},
	journal = {IEEE Robotics and Automation Letters},
	author = {Wu, Xuwei and Kirner, Annika and Garofalo, Gianluca and Ott, Christian and Kotyczka, Paul and Dietrich, Alexander},
	month = oct,
	year = {2022},
	pages = {12307--12314},
}

@inproceedings{traversaro_identification_2016,
	title = {Identification of fully physical consistent inertial parameters using optimization on manifolds},
	doi = {10.1109/IROS.2016.7759801},
	booktitle = {{IEEE}/{RSJ} {International} {Conference} on {Intelligent} {Robots} and {Systems}},
	author = {Traversaro, Silvio and Brossette, Stanislas and Escande, Adrien and Nori, Francesco},
	month = oct,
	year = {2016},
	pages = {5446--5451},
}

@article{pagilla_adaptive_2000,
	title = {Adaptive control of time-varying mechanical systems: analysis and experiments},
	volume = {5},
	doi = {10.1109/3516.891052},
	number = {4},
	journal = {IEEE/ASME Transactions on Mechatronics},
	author = {Pagilla, P.R. and Yu, B. and Pau, K.L.},
	month = dec,
	year = {2000},
	pages = {410--418},
}

@article{ghorbel_uniform_1998,
	title = {On the uniform boundedness of the inertia matrix of serial robot manipulators},
	volume = {15},
	doi = {10.1002/(SICI)1097-4563(199812)15:1<17::AID-ROB2>3.0.CO;2-V},
	number = {1},
	journal = {Journal of Robotic Systems},
	author = {Ghorbel, Fathi and Srinivasan, B. and Spong, Mark W.},
	year = {1998},
	pages = {17--28},
}

@article{johansson_adaptive_1990,
	title = {Adaptive control of robot manipulator motion},
	volume = {6},
	doi = {10.1109/70.59359},
	number = {4},
	journal = {IEEE Transactions on Robotics and Automation},
	author = {Johansson, R.},
	month = aug,
	year = {1990},
	pages = {483--490},
}

@phdthesis{lee_geometric_2019,
	title = {Geometric {Methods} for {Dynamic} {Model}-{Based} {Identification} and {Control} of {Multibody} {Systems}},
	school = {Seoul National University},
	author = {Lee, Taeyoon},
	month = jul,
	year = {2019},
}

@inproceedings{cho_recursive_2024,
	title = {Recursive {Least} {Squares} with {Log}-{Determinant} {Divergence} {Regularisation} for {Online} {Inertia} {Identification}},
	doi = {10.1109/ICRA57147.2024.10610389},
	urldate = {2025-11-10},
	booktitle = {{IEEE} {International} {Conference} on {Robotics} and {Automation}},
	author = {Cho, Namhoon and Lee, Taeyoon and Shin, Hyo-Sang},
	month = may,
	year = {2024},
	pages = {12578--12584},
}

@inproceedings{kaufmann_NICL_2026,
	booktitle = {23rd {IFAC} World Congress},
	title = {Integral {Concurrent} {Learning} for {Natural} {Adaptive} {Control} of {Robotic} {Manipulators}},
	author = {Kaufmann, Tom and Reger, Johann},
	month = aug,
	year = {2026},
}

@book{folland_real_1999,
	title = {Real {Analysis}: {Modern} {Techniques} and {Their} {Applications}},
	publisher = {John Wiley \& Sons},
	author = {Folland, Gerald B.},
	year = {1999},
}

\end{document}